\documentclass[11pt]{article}
\usepackage[dvipdfm]{graphicx, color}
\usepackage{amssymb, amsmath,amsthm}

\newtheorem{thm}{{\bf Theorem}}[section]
\newtheorem{lem}[thm]{{\bf Lemma}}
\newtheorem{cor}[thm]{{\bf Corollary}}
\newtheorem{prop}[thm]{{\bf Proposition}}

\newtheorem{claim}[thm]{Claim} 
\newtheorem{rem}[thm]{Remark}
\newtheorem{ex}[thm]{Example}
\newtheorem{prob}[thm]{Problem}

\newtheorem{conj}[thm]{Conjecture}
\newtheorem{definition}[thm]{Definition}
\numberwithin{equation}{section}

\begin{document}
\title{Pseudo-Anosov braids with small entropy and the magic $3$-manifold}
 \author{Eiko Kin\footnote{The author is partially supported by Grant-in-Aid for Young Scientists (B) (No. 20740031), 
MEXT, Japan}\hspace{1mm} 
and Mitsuhiko Takasawa}
\maketitle

\begin{abstract}
 We consider a surface bundle over the circle, the so called magic manifold  $M$. 
 We determine homology classes whose minimal representatives are genus $0$ fiber surfaces for $M$, and describe their  monodromies by braids. 
Among those  classes whose representatives have $n$ punctures for each $n$, we decide which one realizes the minimal entropy. 
We show that for each $n \ge 9$ (resp. $n= 3,4,5,7,8$), 
there exists a pseudo-Anosov homeomorphism $\Phi_n: D_n \rightarrow D_n$ 
with the smallest known entropy (resp. the smallest entropy) 
which occurs as the monodromy on an $n$-punctured disk fiber for the Dehn filling of $M$. 
A pseudo-Anosov homeomorphism $\Phi_6: D_6 \rightarrow D_6$ 
with the smallest  entropy occurs as the monodromy on a $6$-punctured disk fiber for $M$. 
\end{abstract}

\medskip


\section{Introduction}

Let  $\mathcal{M}(\Sigma)$ be the mapping class group of an orientable surface $\Sigma=\Sigma_{g,p}$ of genus $g$ with $p$ punctures. 
Assuming that $3g-3+p \ge 1$, 
elements of $\mathcal{M}(\Sigma)$ are classified into three types: periodic, pseudo-Anosov  and reducible \cite{Thurston2}. 
There exist  two  numerical invariants of  pseudo-Anosov mapping classes $\phi$. 
One  is the entropy $\mathrm{ent}(\phi)$ which is the logarithm of the dilatation $\lambda(\phi)$. 
The other  is the volume $\mathrm{vol}(\phi)$ which comes from
the hyperbolization theorem by Thurston \cite{Thurston3}.
His theorem asserts that  $\phi$ is pseudo-Anosov if and only if its mapping torus 
$${\Bbb T}(\phi) =\Sigma \times [0,1]/ \sim ,$$ 
 where $\sim$ identifies $(x,1)$ with $(f(x),0)$ for any representative $f \in \phi$, is hyperbolic. 
 We denote   the volume of ${\Bbb T}(\phi)$ by $\mathrm{vol}(\phi)$. 
 
 Let $\mathcal{M}^{\mathrm{pA}}(\Sigma)$ be the set of pseudo-Anosov elements of $\mathcal{M}(\Sigma)$. 
 Fixing  $\Sigma $, the dilatation $\lambda(\phi)$ for $\phi \in \mathcal{M}^{\mathrm{pA}}(\Sigma)$ is known to be an algebraic integer with a 
 bounded degree depending only on $\Sigma$. 
 The set of dilatations $\lambda(\phi)$ for $\phi \in \mathcal{M}^{\mathrm{pA}}(\Sigma)$ bounded by each constant from above is finite, see \cite{Ivanov}. 
 In particular the set 
$$\mathrm{Dil}(\Sigma) = \{\lambda(\phi)>1\ |\ \phi \in \mathcal{M}^{\mathrm{pA}}(\Sigma)\}$$ 
achieves its infimum  $\lambda(\Sigma)$. 

We turn to  volume. 
The set 
$$\{v \ |\ v\ \mbox{is\ the\ volume\ of\ a\ hyperbolic\ }3\mbox{-manifold}\}$$ 
is a well-ordered closed subset  of ${\Bbb R}$ of order type $\omega^{\omega}$ \cite{Thurston}. 
In particular any subset  achieves its infimum. 
Let $ \mathrm{vol}(\Sigma) = \min\{\mathrm{vol}(\phi)\ |\ \phi \in \mathcal{M}^{\mathrm{pA}}(\Sigma)\}$. 
It is of interest to compute  $\lambda(\Sigma)$ (resp.  $\mathrm{vol}(\Sigma)$) 
and to determine the mapping class realizing the minimum. 
Another problem related to the minimal dilatation (resp. minimal volume) are as follows. 
For a non-negative integer $c$, we set 
\begin{eqnarray*}
\lambda(\Sigma;c) &=& \min\{\lambda(\phi) \ |\  \phi \in \mathcal{M}^{\mathrm pA}(\Sigma),\ {\Bbb T}(\phi)\ \mbox{has\ }c\mbox{\ cusps}\}, 
\\
\mathrm{vol}(\Sigma;c)&=& \min\{\mathrm{vol}(\phi)\ |\  \phi \in \mathcal{M}^{\mathrm{pA}}(\Sigma),\  {\Bbb T}(\phi)\ \mbox{has\ }c\mbox{\ cusps}\}.
\end{eqnarray*}
A problem is to compute $\lambda(\Sigma;c) $ (resp. $\mathrm{vol}(\Sigma;c)$) and to find a mapping class realizing the minimum.

\begin{figure}
\begin{center}
\includegraphics[width=1in]{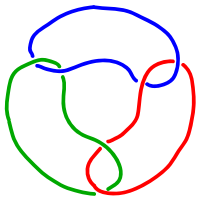}
\caption{$3$ chain link $\mathcal{C}_3$.}
\label{fig_3chain}
\end{center}
\end{figure}

In \cite{KKT}, the authors and S.~Kojima obtain experimental results concerning the entropy and volume. 
In the case  the mapping class group $\mathcal{M}(D_n)$ of an $n$-punctured disk $D_n$, they observe that 
for many pairs $(n,c)$, there exists  a mapping class simultaneously reaching both $\lambda(D_n;c)$ and $\mathrm{vol}(D_n;c)$. 
Experiments tell us that in case $c=3$, the mapping tori reaching both minima are homeomorphic to 
the magic manifold $M_{\mathrm{magic}}$ 
which is  the exterior of the $3$ chain link $\mathcal{C}_3$  illustrated in Figure~\ref{fig_3chain}. 
Moreover when $c=2$, it is observed that  there exists a mapping class $\phi$ realizing both  $\lambda(D_n;2)$ and $\mathrm{vol}(D_n;2)$  and 
its mapping torus ${\Bbb T}(\phi)$ is homeomorphic to a Dehn filling of  $M_{\mathrm{magic}}$  along one cusp. 
This study motivates the present paper which concerns the fibrations in $M_{\mathrm{magic}}$. 
The magic manifold  has the smallest known volume among orientable hyperbolic $3$-manifolds having $3$ cusps. 
Many manifolds having at most $2$ cusps with small volume are obtained from $M_{\mathrm{magic}}$ by Dehn fillings, see \cite{MP}.  
Also,  some  important examples for the study of the exceptional Dehn fillings can be obtained from the Dehn fillings of  $M_{\mathrm{magic}}$, see \cite{Gordon}.

Let $M$ be a hyperbolic $3$-manifold with boundary  which fibers over the circle. 
We assume that $M$ admits infinitely many different fibrations. 
Thurston introduced  the norm function $X_T: H_2(M, \partial M; {\Bbb R}) \rightarrow {\Bbb R}$, and 
showed that the unit ball $U$ with respect to $X_T$ is a compact, convex polyhedron \cite{Thurston1}. 
He described the relation between the function $X_T$ and fibrations of $M$ as follows. 
For each fiber $F$ of $M$, the homology class $[F] \in H_2(M, \partial M; {\Bbb R})$ lies in the open cone $int(C_{\Delta})$ with the origin over a top 
dimensional face $\Delta$ of $\partial U$. 
Conversely for any integral class $a \in int(C_{\Delta})$, there exists a fiber $F$ of $M$ representing $a$.  
Using this description of the fibers, the entropy function $ \mathrm{ent}(\cdot): int(C_{\Delta}) \rightarrow {\Bbb R}$ can be defined as follows. 
For each primitive integral class $a \in int(C_{\Delta})$,  the monodromy  $\Phi_a: F_a \rightarrow F_a$ on a connected surface $F_a$ 
representing $a$ is  pseudo-Anosov, and 
one defines the entropy  of $a$ by $\mathrm{ent}(a) = \log (\lambda (\Phi_a))$. 
Fried proves that this function defined  on primitive integral classes admits a unique continuous extension to a homogeneous function on $int(C_{\Delta})$ 
\cite{Fried}.

One sees that $\mathcal{M}(D_n)$  is isomorphic to the subgroup of 
 $\mathcal{M}(\Sigma_{0,n+1})$ consisting of the elements   which  fix the  puncture of $\Sigma_{0,n+1}$. 
 By using the natural surjective homomorphism $\Gamma: B_n \rightarrow \mathcal{M}(D_n) $ from the $n$-braid group $B_n$ to $\mathcal{M}(D_n)$, 
 one represents each element of $\mathcal{M}(D_n)$ by an $n$-braid.  
 A braid $b$ is called {\it pseudo-Anosov} if $\Gamma(b)$ is a pseudo-Anosov mapping class. 
 If this is the case, the dilatation $\lambda(b)$ of $b$  is defined by the dilatation of $\Gamma(b)$. 
 Let $T_{m,p}$ be the following $m$-braid for $m \ge 3$ and $p \ge 1$ which is a main example in this paper. 
$$T_{m,p} = (\sigma_1^2  \sigma_2 \sigma_3 \cdots \sigma_{m-1})^p \sigma_{m-1}^{-2}.$$
For example, $T_{6,2} = (\sigma_1^2 \sigma_2 \sigma_3 \sigma_4 \sigma_5)^2 \sigma_{5}^{-2} = 
\sigma_1^2 \sigma_2 \sigma_3 \sigma_4 \sigma_5 \sigma_1^2 \sigma_2 \sigma_3 \sigma_4 \sigma_5^{-1}$ 
(Figure~\ref{fig_t_reducible}(left)).  
The braid $T_{m,p}$ is a  horseshoe braid if $\gcd(p,m-1) =1$ and $1<p \le \frac{m-1}{2}$ (Proposition~\ref{prop_horseshoe-type}). 
If $\gcd(m-1,p)=1$, then the mapping torus ${\Bbb T}(\Gamma(T_{m,p}))$ is homeomorphic to $M_{\mathrm{magic}}$ (Corollary~\ref{cor_tmp-monodromy}). 
Otherwise $\Gamma(T_{m,p})$ is  reducible. 
 We set 
$$\mathcal{M}_{\mathrm{magic}}^{n}= 
\{\phi \in \mathcal{M} (\Sigma_{0,n})\ |\ {\Bbb T}(\phi) \ \mbox{is\ homeomorphic\  to\ } M_{\mathrm{magic}} \}.$$
Let us define an integral polynomial 
$$f_{(x,y,z)}(t)=  t^{x+y-z}-t^x - t^y - t^{x-z}- t^{y-z}+1.$$
The following, our main theorem,  states that for all $n$ but $6$ and $8$, 
 the minimum among the dilatations of  $ \phi \in \mathcal{M}_{\mathrm{magic}}^{n}$ is realized by  $\Gamma(T_{n-1,p})$ 
 for some $p=p(n)$ and it is computed as the largest real root of one of the polynomials $f_{(x,y,z)}(t)$.  

\begin{thm} 
\label{thm_main}
For each $n \ge 4$, the minimum among the dilatations of  $\phi \in \mathcal{M}_{\mathrm{magic}}^{n}$ is realized by: 
\begin{description}
\item[(1)] 
$\Gamma(T_{2k,2}) $ in case  $n=2k+1$ for $k \ge 2$.
The dilatation $\lambda(T_{2k,2})$ equals the largest real root of 
$$f_{(k-1,k,0)}(t)= t^{2k-1}- 2(t^{k-1}+ t^k)+1.$$
\item[(2-i)] 
$\Gamma(\sigma_1 \sigma_2^2 \sigma_3 \sigma_4) $ in case $n=6$. 
The dilatation 
$\lambda(\sigma_1 \sigma_2^2 \sigma_3 \sigma_4) \approx 2.08102$ equals  the largest real root of 
$$f_{(3,2,1)}(t)=   t^4-t^3-2t^3-t+1.$$
\item[(2-ii)] 
$\Gamma(T_{4k+1,2k-1}) $ in case  $n=4k+2$ for $k \ge 2$. 
The dilatation $\lambda(T_{4k+1,2k-1})$ equals the largest real root of 
$$f_{(2k+1, 2k-1,0)}(t)=    t^{4k}-2(t^{2k-1}+ t^{2k+1}) +1.$$
\item[(3a-i)] 
$\Gamma(T_{3,1}) $ in case $n=4$. 
The dilatation  $\lambda(T_{3,1}) \approx 3.73205$ equals the largest real root of 
$$f_{(1,1,0)}(t)=    t^2-4t+1.$$
\item[(3a-ii)] 
$\Gamma(T_{8k+3,2k+1})$ in case  $n= 8k+4$ for $k \ge 1$. 
The dilatation $\lambda(T_{8k+3,2k+1})$ equals the largest real root of 
$$f_{(4k-1, 4k+3,0)}(t)=    t^{8k+2}-2(t^{4k-1}+ t^{4k+3})+1.$$
\item[(3b-i)] 
$\Gamma(b) $ in case $n=8$, 
where 
$$b= \sigma_1^{-1} \sigma_2^{-1} \sigma_3^{-1} \sigma_4^{-1} \sigma_5^{-1} 
\sigma_6^{-1} \sigma_1^{-1} \sigma_2^{-1} \sigma_3^{-1} \sigma_4^{-1} \sigma_1^{-1} \sigma_2^{-1} \sigma_3^{-1} \in B_7.$$
The dilatation $\lambda(b) \approx 1.72208$ equals the largest real root of 
$$f_{(5,3,2)}(t)=   t^6-t^5-2t^3-t+1.$$
\item[(3b-ii)] 
$\Gamma(T_{8k+7,2k+1}) $ in case  $n= 8(k+1)$ for $k \ge 1$.
The dilatation $\lambda(T_{8k+7,2k+1})$ equals the largest real root of 
$$f_{(4k+5, 4k+1,0)}(t)=   t^{8k+6}-2(t^{4k+1}+ t^{4k+5})+1.$$
\end{description} 
Moreover the above mapping class realizing the minimal dilatation among elements of $\mathcal{M}_{\mathrm{magic}}^{n}$ is unique up to conjugacy. 
\end{thm}

Forgetting the 1st strand of  $T_{m,p}$, one obtains the $(m-1)$-braid, call it $T'_{m,p}$.  
For example, $T_{6,2}'= \sigma_1 \sigma_2 \sigma_3 \sigma_4 \sigma_1 \sigma_2 \sigma_3 \sigma_4^{-1}$.  
The families of braids $\{T'_{m,p}\}$ and $\{T_{m,p}\}$ contain examples of $\ell$ strands ($\ell= 3,4, \cdots, 8$) with the smallest dilatation. 
(See Section~\ref{subsection_pA}.) 
Hironaka-Kin (resp. Venzke) found candidates with the smallest dilatation $\lambda(D_n)$  for $n$ odd (resp. $n$ even), see \cite{HK} (resp. \cite{Venzke}).   
All the braids  in Theorem~\ref{thm_main}(1)(2-ii)(3a-ii)(3b-ii) relate to those examples. 
The braid $T_{2k,2}'$ (with odd strands) is conjugate to the braid $\sigma_{(k)}$   with the smallest known dilatation 
found by Hironaka-Kin.  (See Theorem~\ref{thm_main}(1).)  
For the braid $T'_{m,p}$ (with even strands) obtained from $T_{m,p}$ in (2-ii), (3a-ii) or (3b-ii) of Theorem~\ref{thm_main}, 
the mapping class $\Gamma(T'_{m,p})$ is conjugate to the one given by Venzke. 
(See Section~\ref{subsection_pA}.)

Work of Farb-Leininger-Margalit \cite{FLM} together with a result in \cite{HK} implies  that 
there exists a complete, noncompact, finite volume, hyperbolic $3$-manifold $M'$ 
with the following property: 
there exist Dehn fillings of $M'$ giving an infinite sequence of fiberings over $S^1$, with fibers $D_{n_i}$ having  $n_i$ punctures 
with  $n_i \to \infty$, and with the monodromy $\Phi_i: D_{n_i} \rightarrow D_{n_i}$ so that $\lambda(D_{n_i}) = \lambda(\Phi_i)$.  
The magic manifold is a potential example which could satisfy this property.

In  \cite{AD,Hironaka,KT1}, one can find 
pseudo-Anosovs on closed surfaces $\Sigma_g$ of genus $g$ with small dilatation which occur as monodromies on fibers for Dehn fillings of $M_{\mathrm{magic}} $. 
Using those pseudo-Anosovs,  
Hironaka \cite{Hironaka}, Aaber-Dunfield \cite{AD}, and the authors \cite{KT1} independently  proved that 
$$\displaystyle \lim_{g \to \infty} \sup g \log  \lambda(\Sigma_g) \le \log(\tfrac{3+ \sqrt{5}}{2}).$$

This paper is organized as follows. 
Section~\ref{section_notation} reviews basic facts. 
Section~\ref{section_magic-manifold} contains the proof of Theorem~\ref{thm_main}. 
For the proof, we first compute the Teichm\"{u}ler polynomial, introduced by McMullen \cite{McMullen}, which determines the 
entropy function for $M_{\mathrm{magic}}$ (Theorem~\ref{thm_Fried-Oertel-Poly}). 
Then we find all the homology classes whose representatives are genus $0$ fiber surfaces (Corollary~\ref{cor_genus0}). 
We study the asymptotic behaviors of the normalized entropy function 
$\overline{\mathrm{ent}}(\cdot) = X_T(\cdot) \mathrm{ent}(\cdot)$ (Theorem~\ref{thm_asym-P}). 
This tells us  which class realizes the minimal dilatation 
among homology classes whose representatives are genus $0$ fiber surfaces with $n$ punctures 
(Proposition~\ref{prop_minimal-homology}). 
We finally describe the monodromies for these fiber surfaces by using braids (Propositions~\ref{prop_mini-monodromy} and \ref{prop_two-monodromies}). 
In Section~\ref{section_further},  we discuss pseudo-Anosov braids with small dilatation. 
We also find a relation between the horseshoe map and the braids $T_{m,p}$ (Proposition~\ref{prop_horseshoe-type}). 
\medskip

\noindent
Acknowledgments: 
We would like to thank Shigeki Akiyama and Takuya Sakasai who showed us the proof of  Proposition~\ref{prop_genus0}. 
We would like to thank  Sadayoshi Kojima and Makoto Sakuma for a great deal of encouragement. 
We also would like to thank the referees for valuable comments and suggestions.

\section{Notation and basic facts}
\label{section_notation}

\subsection{Mapping class group}

The {\it mapping class group} $\mathcal{M}(\Sigma)$ is the group of isotopy classes of orientation preserving homeomorphisms of $\Sigma$, 
where the group operation is induced by composition of homeomorphisms. 
An element of the mapping class group  is called a {\it mapping class}. 

A homeomorphism $\Phi: \Sigma \rightarrow \Sigma$ is {\it pseudo-Anosov}  
 if  there exists a constant $\lambda= \lambda(\Phi)>1$ called the {\it dilatation of} $\Phi$  
 and there exists a pair of transverse measured foliations $\mathcal{F}^s$ and $\mathcal{F}^u$ such that 
 $$\Phi(\mathcal{F}^s)= \tfrac{1}{\lambda} \mathcal{F}^s \ \mbox{and}\  \Phi(\mathcal{F}^u)= \lambda \mathcal{F}^u.$$ 
 In this case the mapping class  $\phi=[\Phi]$ is called pseudo-Anosov. 
 We define the dilatation of $\phi$, denoted by $\lambda(\phi)$, to be the dilatation of $\Phi$.

 The (topological)  entropy $\mathrm{ent}(f)$ is a measure of the complexity of a continuous self-map $f$ on a compact manifold, see  \cite{Walters}. 
 For a pseudo-Anosov homeomorphism $\Phi: \Sigma \rightarrow \Sigma$, 
the equality 
$$\mathrm{ent}(\Phi)= \log (\lambda(\Phi))$$ 
holds  and $\mathrm{ent}(\Phi)$ attains the minimal entropy among all  homeomorphisms  
which are isotopic to $\Phi$, see \cite[Expos\'{e} 10]{FLP}.  
We denote by  $\mathrm{ent}(\phi)$  this characteristic number.  
Using the Euler characteristic $\chi(\Sigma)$, 
we define the {\it normalized entropy} and {\it normalized dilatation}  of $\phi$ by 
$\overline{\mathrm{ent}}(\phi)= |\chi(\Sigma)| \mathrm{ent}(\phi)$ and 
$\overline{\lambda}(\phi)= |\chi(\Sigma)| \lambda(\phi)$.

\begin{figure}[htbp]
\begin{center}
\includegraphics[width=3.7in]{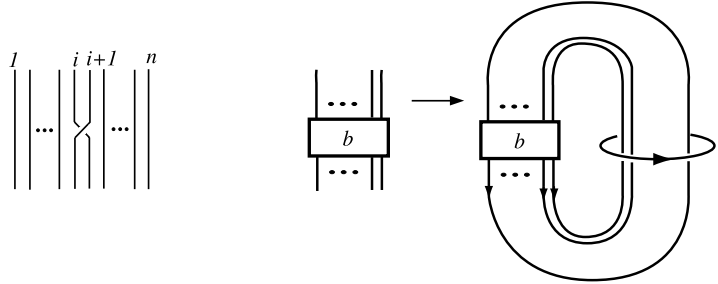}
\caption{(left) generator $\sigma_i$. (right) braid $b \rightarrow $ braided link $\overline{b}$.}
\label{fig_artin}
\end{center}
\end{figure}

We recall the surjective homomorphism 
$$\Gamma: B_n \rightarrow \mathcal{M}(D_n)$$
which sends the  Artin generator $\sigma_i$ for $i \in \{1, \cdots, n-1\}$ (see Figure~\ref{fig_artin}(left)) to $\hat{t}_i$, 
where $\hat{t}_i$ is the mapping class which represents the positive half twist about the arc from the $i$th puncture to the $(i+1)$st puncture.  
The kernel of $\Gamma$ is the center of $B_n$ which is generated by the full twist 
$(\sigma_1 \sigma_2 \cdots \sigma_{n-1})^n$. 
By replacing the boundary of $D_n$ with the $(n+1)$st  puncture, 
the injective homomorphism from 
$\mathcal{M}(D_n)$ to  $\mathcal{M}(\Sigma_{0,n+1})$ is induced. 
In the rest of the paper we regard an element of $\mathcal{M}(D_n)$  as an element of  $\mathcal{M}(\Sigma_{0,n+1})$.

We say that a braid $b \in B_n$ is {\it pseudo-Anosov} if $\Gamma(b) \in \mathcal{M}(D_n)$ is pseudo-Anosov. 
Is the  case,  
$\mathrm{vol}(\Gamma(b))$ equals the hyperbolic volume of  the exterior  of the link  $ \overline{b}$ in $S^3$, where 
$\overline{b}$ is a union of the closed braid of $b$ and the braid axis. 
Our convention of the orientation of $\overline{b}$ is given by Figure~\ref{fig_artin}(right).

\subsection{Roots of polynomials}

Let $f(t)$ be an integral polynomial of degree $d$. 
The reciprocal of $f(t) $ is  $f_*(t) = t^d f(1/t)$. 
We denote by $\lambda(f)$, the maximal absolute value of the roots of $f(t)$.

Let $R(t)$ be a monic integral polynomial and let $S(t)$  be  an integral polynomial. 
We set 
$$Q_{n,\pm}(t)= t^n  R(t) \pm S(t)$$ 
for each integer $n \ge 1$. 
In case  $S(t)= R_*(t)$, we call $Q_{n,\pm}(t)= t^n  R(t) \pm R_*(t)$  the {\it Salem-Boyd polynomial} {\it associated to} $R(t)$.

\begin{lem}
\label{lem_asymptotic-root1} 
Let $Q_{n,\pm}(t)= t^n  R(t) \pm S(t)$. 
Suppose that $R(t)$ has a root outside the unit circle.  
Then, the roots of $Q_{n,\pm}(t)$ outside the unit circle converge to those of $R(t)$ counting multiplicity as $n$ goes to $\infty$.  
In particular, 
$$ \lambda(R) = \lim_{n \to \infty} \lambda(Q_{n,\pm}).$$
\end{lem}

\noindent
The proof of Lemma~\ref{lem_asymptotic-root1} can be found in \cite[Lemma~2.5]{KT}. 

\subsection{Hyperbolic surface bundle over the circle}

Let $M$ be an irreducible, atoroidal and oriented $3$-manifold with boundary $\partial M$ (possibly $\partial M = \emptyset$). 
Thurston discovered a norm function $X_T: H_2(M, \partial M; {\Bbb R}) \rightarrow {\Bbb R}$ (see  \cite{Thurston1}).  
In case  $M$ is a surface bundle over the circle, he described a relation between  $X_T$ and fibrations of $M$  
which we record as Theorem~\ref{thm_norm-fibration} below.

\subsubsection{Thurston norm} 

The norm function 
$X_T: H_2(M, \partial M; {\Bbb R}) \rightarrow {\Bbb R}$ has the property that 
for any integral class $a  \in H_2(M, \partial M; {\Bbb R})$, 
$$X_T(a)= \min_F \{- \chi(F)\},$$  
where 
 the minimum is taken over all  oriented surface $F$ embedded in $M$, satisfying $a= [F]$,  
 with no components  of non-negative Euler characteristic. 
The surface $F$ which realizes this minimum  is called  a {\it minimal representative} of $a $.  
For a rational class $a \in H_2(M, \partial M; {\Bbb R})$, take a rational number $r $ so that 
$ra$ is an integral class. 
Then $X_T(a)$ is defined to be 
$$X_T(a) = \tfrac{1}{|r|} X_T(ra).$$ 
The function $X_T$ defined on rational classes  admits a unique continuous extension to 
$H_2(M, \partial M; {\Bbb R}) $ which is linear on the ray though the origin. 
The unit ball $U= \{a \in H_2(M, \partial M; {\Bbb R})\ |\ X_T(a) \le 1\}$ is a compact, convex polyhedron \cite{Thurston1}.

The following notations are needed to describe how fibrations of $M$ are related to the Thurston norm. 

\begin{itemize}
\item A top dimensional face in the boundary $\partial U$ of the unit ball $U$ is denoted by $\Delta$, and its open face is denoted by $int(\Delta)$. 
\item The open cone with the origin  over $\Delta$ is denoted by $int(C_{\Delta})$. 
\item The set of integral classes of $int(C_{\Delta})$ is denoted by $int(C_{\Delta}(\Bbb {Z})) $, and 
the set of rational classes of   $int(C_{\Delta})$ is denoted by $int(C_{\Delta}(\Bbb {Q})) $.  
\end{itemize}

\begin{thm}[\cite{Thurston1}] 
\label{thm_norm-fibration}
Suppose that $M$ is a surface bundle over the circle and let $F$ be a fiber. 
Then there exists a top dimensional  face $\Delta$ satisfying the following. 
\begin{description}
\item[(1)] $[F] \in int(C_{\Delta}(\Bbb Z))$. 
\item[(2)]  
For any $a \in int(C_{\Delta}(\Bbb Z))$, 
a minimal  representative $E$ of $a$ is a  fiber of fibrations of $M$.  
\end{description}
\end{thm}
\noindent
The face $\Delta$  in Theorem~\ref{thm_norm-fibration} is called the {\it fiber face}. 
For the fiber face $\Delta$, it follows that $a \in int(C_{\Delta}(\Bbb Z))$ is a primitive integral class if and only if 
a minimal  representative $E$ of $a$ is  connected.

It is known that if $a_0 \in H_2(M, \partial M; {\Bbb Z})$ has a representative $F$ which is a fiber of the fibration of $M$, 
then any incompressible surface which  represents  $a_0$ is isotopic to the fiber $F$, see \cite{Thurston1}. 
In particular $F$ is a minimal representative of $a_0$.  
Thus a minimal representative of $a_0$ is unique up to isotopy. 

\subsubsection{Entropy function}

Suppose that $M$ is a hyperbolic  surface bundle over the circle. 
We fix a fiber face $\Delta$ for $M$. 
The {\it entropy function} 
$ \mathrm{ent}(\cdot): int(C_{\Delta}({\Bbb Q})) \rightarrow {\Bbb R}$ introduced by Fried in  \cite{Fried} is defined as follows. 
The minimal  representative $F_a$ for $a \in int(C_{\Delta}({\Bbb Z}))$  is a fiber of fibrations of $M$. 
Let  $\Phi_a : F_a \rightarrow F_a$ be the monodromy. 
Since $M$ is a hyperbolic manifold, the mapping class $\phi_a = [\Phi_a]$  must be  pseudo-Anosov. 
The entropy  $\mathrm{ent}(a) $ and  dilatation $\lambda(a)$ are defined as the entropy and dilatation of $\phi_a$, respectively.
For a rational number $r $ and an integral class $a$, 
the entropy $\mathrm{ent}(ra) $ is defined by $ \frac{1} {|r|}  \mathrm{ent}(a) .$
Notice that $X_T(\cdot) \mathrm{ent}(\cdot):  int(C_{\Delta}({\Bbb Q})) \rightarrow {\Bbb R}$ 
is constant on each ray through the origin.  
We call $X_T(a) \mathrm{ent}(a)$ and $X_T(a) \lambda(a)$ the  {\it normalized entropy} and {\it normalized dilatation} of $a$.

We recall an important property of  the entropy function proved by Matsumoto and independently McMullen.

\begin{thm}[\cite{Matsumoto,McMullen}] 
\label{thm_strictly-concave}
The function 
$\tfrac{1}{\mathrm{ent}(\cdot)}: int(C_{\Delta}({\Bbb Q})) \rightarrow {\Bbb R}$ 
is strictly concave. 
\end{thm}

\noindent
By Theorem~\ref{thm_strictly-concave}, the function  $\mathrm{ent}(\cdot)$ on $int(C_{\Delta}({\Bbb Q}))$  admits a unique continuous extension to 
$\mathrm{ent}(\cdot): int(C_{\Delta}) \rightarrow {\Bbb R}$. 

Since $\mathrm{ent}(a)$ goes to $ \infty$ as $a $ goes to a point on the boundary $\partial \Delta$ (see \cite{Fried}),  
Theorem~\ref{thm_strictly-concave} implies the normalized entropy function
$$ \overline{\mathrm{ent}}(\cdot)= X_T(\cdot) \mathrm{ent}(\cdot):  int(C_{\Delta}) \rightarrow {\Bbb R}$$ 
has the minimum at a unique ray through the origin. 
In other words $\mathrm{ent}(\cdot)$ has the minimum at a unique point  of  $int(\Delta)$.  
The following question was posed by McMullen \cite[p~542]{McMullen}.

\begin{prob}
\label{prob_minimal-ray}
On which ray does it attain the minimal normalized entropy with respect to the fiber face?
Is the minimum always attained at a rational class of $int(\Delta)$?
\end{prob}

\noindent 
We solve this problem for the magic manifold in  Section~\ref{subsection_teich}.

\section{Magic manifold}
\label{section_magic-manifold}

\subsection{Fiber face}
\label{subsection_FiberFace}

\begin{figure}
\begin{center}
\includegraphics[width=4.3in]{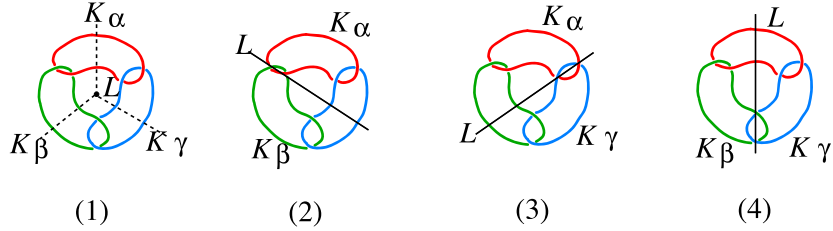}
\caption{axis $\mathtt{L}$ for periodic map}
\label{fig_symmetry}
\end{center}
\end{figure}

\begin{figure}[htbp]
\begin{center}
\includegraphics[width=4in]{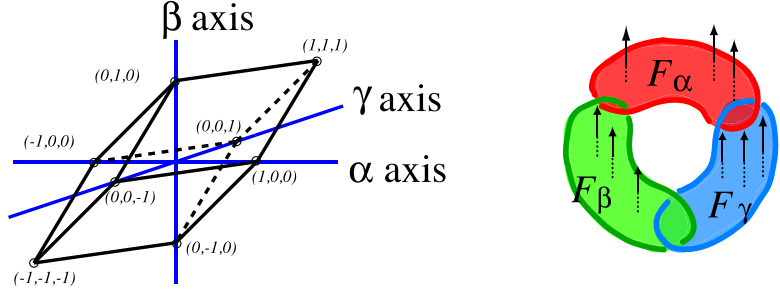}
\caption{(left) Thurston unit ball. (right) $F_{\alpha}$, $F_{\beta}$, $F_{\gamma}$. 
(arrows indicate the normal direction of oriented surfaces.)}
\label{fig_poly}
\end{center}
\end{figure}

Let $\mathtt{L}$ be one of the four lines  in $S^3$ depicted  in Figure~\ref{fig_symmetry}(1),(2),(3) and (4). 
For each $\mathtt{L}$, there exist an integer $n$ and  a periodic map $f: (S^3, {\mathcal C}_3) \rightarrow (S^3, {\mathcal C}_3)$ such that 
$f$ is a $2 \pi/n$ rotation with respect to $\mathtt{L}$. 
Such symmetry of ${\mathcal C}_3$ reflects the shape of the Thurston unit ball. 
Let $K_{\alpha}$, $K_{\beta}$ and  $K_{\gamma}$ be the components of $\mathcal{C}_3$ such that 
$K_{\alpha}$ (resp. $K_{\beta}$, $K_{\gamma}$) bounds the oriented twice-punctured disk $F_{\alpha}$ (resp. $F_{\beta}$, $F_{\gamma}$) 
in $M_{\mathrm{magic}}$ whose normal direction is indicated as in Figure~\ref{fig_poly}(right). 
Those oriented surfaces induce the orientation of $\mathcal{C}_3$. 
Let $\alpha = [F_{\alpha}]$, $ \beta = [F_{\beta}]$, and $\gamma= [F_{\gamma}]$. 
In \cite{Thurston1}, Thurston computes the unit ball $U$ which is the the parallelepiped with vertices 
$\pm \alpha = (\pm 1,0,0)$, $\pm \beta = (0,\pm 1,0)$, $\pm \gamma= (0, 0, \pm1)$, $\pm(\alpha + \beta + \gamma)$, 
see Figure~\ref{fig_poly}(left). 
The set $\{\alpha, \beta, \gamma\}$ is a basis of $H_2(M_{\mathrm{magic}}, \partial M_{\mathrm{magic}}; {\Bbb Z})$. 

\begin{figure}[htbp]
\begin{center}
\includegraphics[width=5in]{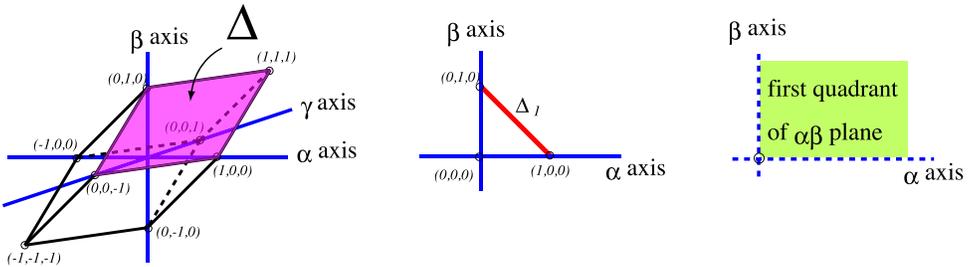}
\caption{(left) fiber face $\Delta$. (center) $\Delta_1 \subset \Delta$. (right) $C_{\Delta_1} \subset int(C_{\Delta})$. }
\label{fig_face}
\end{center}
\end{figure}

The magic manifold is  a surface bundle over the circle as we will see it later.  
The symmetry of ${\mathcal C}_3$ tells us  that every top dimensional face is a fiber face. 
We (arbitrarily)  pick the shaded fiber face $\Delta$ as in Figure~\ref{fig_face}(left) 
with vertices $\alpha = (1,0,0)$, $\alpha+ \beta + \gamma= (1,1,1)$, $\beta=(0,1,0)$ and $-\gamma= (0,0,-1)$.  
The open face $ int(\Delta) $ is written by 
\begin{equation}
\label{equation_int-delta}
 int(\Delta) = \{x \alpha + y \beta + z \gamma\ |\ x+y-z =1, \ x >0,\  y>0,\   x >z,\   y>z\}.
 \end{equation}
 For  $x \alpha + y \beta + z \gamma \in int(C_{\Delta})$ (not necessarily primitive), 
$$X_T(x \alpha + y \beta + z \gamma) = x + y -z.$$

Let $\mathcal{N}(L)$ be the  regular neighborhood of a link $L$ in $S^3$. 
We denote the tori $ \partial \mathcal{N}(K_{\alpha})$, $ \partial \mathcal{N}(K_{\beta})$, $ \partial \mathcal{N}(K_{\gamma})$ by 
$T_{\alpha}$, $T_{\beta}$, $T_{\gamma}$ respectively. 
Let  $ x \alpha + y \beta + z \gamma $ be a primitive integral class in $int(C_{\Delta})$. 
We denote by $F_{x \alpha + y \beta + z \gamma}$ or $F_{(x,y,z)}$, the minimal representative of $ x \alpha + y \beta + z \gamma $. 
Let us set  
$$\partial_{\alpha} F_{(x,y,z)} = \partial F_{(x,y,z)} \cap T_{\alpha}$$ 
which consists of  the parallel simple closed curves on $T_{\alpha}$. 
We define the subsets $\partial_{\beta} F_{(x,y,z)} $, $\partial_{\gamma} F_{(x,y,z)} \subset  \partial F_{(x,y,z)} $ in the same manner. 
We denote by $\Phi_{(x,y,z)}: F_{(x,y,z)} \rightarrow F_{(x,y,z)}$, the monodromy  on a fiber $F_{(x,y,z)}$. 
It is clear that $\Phi_{(x,y,z)}$ permutes elements of each of the sets $\partial_{\alpha} F_{(x,y,z)}$, 
$\partial_{\beta} F_{(x,y,z)}$ and $\partial_{\gamma} F_{(x,y,z)}$ cyclically. 
Let $\mathcal{F}_{(x,y,z)}$ be the stable foliation for the pseudo-Anosov  $\Phi_{(x,y,z)}$.

\begin{lem}
\label{lem_topological-type}
Let  $ x \alpha + y \beta + z \gamma \in int(C_{\Delta})$ be a primitive integral class. 
The number of the boundary components $\sharp(\partial F_{(x,y,z)})$ is equal to  
the sum of the three greatest common divisors 
$$\gcd(x,y+z)+ \gcd(y,z+x)+ \gcd(z,x+y),$$ 
where $\gcd(0,w)$ is defined by $|w|$.  
More precisely 
\begin{description}
\item[(1)] 
$ \sharp( \partial_{\alpha} F_{(x,y,z)}) = \gcd(x,y+z)$, 
\item[(2)]
$\sharp (\partial_{\beta} F_{(x,y,z)} )=  \gcd(y,z+x)$, 
\item[(3)]
$\sharp( \partial_{\gamma} F_{(x,y,z)}) = \gcd(z,x+y)$. 
\end{description}
\end{lem}

\noindent
{\it Proof.} 
We prove  (1). 
The proof of (2),(3) is similar. 
We have the meridian and longitude basis $\{m_{\alpha}, \ell_{\alpha}\}$ for $T_{\alpha}$. 
Similarly we have the bases $\{m_{\beta}, \ell_{\beta}\}$ for $T_{\beta}$ and 
$\{m_{\gamma}, \ell_{\gamma}\}$ for $T_{\gamma}$. 
We consider the long exact sequence of the homology groups of the pair $(M_{\mathrm{magic}}, \partial M_{\mathrm{magic}})$. 
The boundary map is given by 
\begin{eqnarray*}
\partial_*: H_2(M_{\mathrm{magic}}, \partial M_{\mathrm{magic}}; {\Bbb R})& \rightarrow& H_1(\partial M_{\mathrm{magic}}; {\Bbb R}), 
\\
\alpha &\mapsto& \ell_{\alpha} - m_{\beta}- m_{\gamma}, 
\\
\beta &\mapsto& \ell_{\beta} - m_{\gamma}- m_{\alpha}, 
\\
\gamma &\mapsto& \ell_{\gamma} - m_{\alpha}- m_{\beta}. 
\end{eqnarray*}
Hence 
\begin{equation}
\label{equation_homology}
\partial_*(x \alpha + y \beta + z \gamma) = x \ell_{\alpha} - (y+z) m_{\alpha} + 
y \ell_{\beta} - (z+x) m_{\beta} + z \ell_{\gamma} - (x+y)m_{\gamma}. 
\end{equation}
Since $F_{(x,y,z)}$ is the minimal representative, 
the set $\partial_{\alpha} F_{(x,y,z)}$ is a union of oriented parallel simple closed curves on $T_{\alpha}$ 
whose homology class  equals  $x \ell_{\alpha} - (y+z) m_{\alpha}  \in H_1(T_{\alpha}; {\Bbb R})$, 
see (\ref{equation_homology}). 
Thus the number of the components of $\partial_{\alpha} F_{(x,y,z)}$ equals $\gcd(x,y+z)$. 
$\Box$
\medskip

In Section~\ref{subsection_pA}, we will use the following to see 
$\lambda(T_{m,p}')$ is equal to $\lambda(T_{m,p})$. 

\begin{lem}[\cite{KT1}]
\label{lem_sing-data}
Let $x \alpha+ y \beta + z \gamma$ be as in Lemma~\ref{lem_topological-type}. 
The stable foliation  $\mathcal{F}_{(x,y,z)}$  has 
\begin{description}
\item[(1)]  
$\tfrac{x}{\gcd(x,y+z)}$ prongs  at each element  of $\partial_{\alpha} F_{(x,y,z)}$, 
\item[(2)] 
$\tfrac{y}{\gcd(y,x+z)}$ prongs  at each element  of  $\partial_{\beta} F_{(x,y,z)}$, 
\item[(3)] 
$\tfrac{x+y-2z}{\gcd(z,x+y)}$ prongs  at each element of $\partial_{\gamma} F_{(x,y,z)}$, and 
\item[(4)] 
no singularities in the interior of $F_{(x,y,z)} $. 
\end{description}
\end{lem}

\subsection{Teichm\"{u}ler  polynomial} 
\label{subsection_teich}

We compute the Teichm\"{u}ler  polynomial $P= P_{\Delta}$ with respect to $\Delta$. 
For background, see \cite{McMullen}. 

\begin{figure}[htbp]
\begin{center}
\includegraphics[width=2.7in]{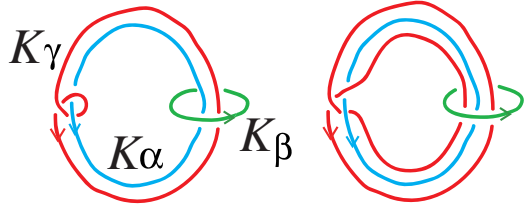}
\caption{(left) $\mathcal{C}_3$. (right) $\overline{b}$. 
(this figure explains how to obtain $H$.)}
\label{fig_chain-homeo}
\end{center}
\end{figure}

A fiber $F= F_{(1,1,0)}$ is homeomorphic to a sphere with $4$ boundary components. 
We now see that 
the monodromy $\Phi= \Phi_{(1,1,0)}$ on $F$ is represented by the $3$-braid $b= \sigma_2 \sigma_1^{-1} \sigma_2$. 
A homeomorphism  $H: S^3 \setminus \mathcal{N}(\mathcal{C}_3) \rightarrow S^3 \setminus \mathcal{N}(\overline{b})$ is given as follows. 
The  link illustrated in Figure~\ref{fig_chain-homeo}(left) is isotopic to $\mathcal{C}_3$. 
We consider the exterior $S^3 \setminus \mathcal{C}_3$ and 
we open the twice-punctured disk $F_{\alpha}$ bounded by $K_{\alpha}$. 
Let $F'_{\alpha}$ and $F''_{\alpha}$ be the resulting twice-punctured
disks obtained from $F_{\alpha}$. 
Reglue $F'_{\alpha}$ and $F''_{\alpha}$ by twisting one of the disks by $360 $ degrees in the clockwise direction. 
Then we obtain the braided link $\overline{b}$  whose exterior $S^3 \setminus \overline{b}$  is homeomorphic to  $S^3 \setminus \mathcal{C}_3$, 
see Figure~\ref{fig_chain-homeo}.

Let $u$ be the meridian of the component of $\overline{b}$ which is the braid axis. 
Let $t_2$ (resp. $t_1$) be the meridian of the component of $\overline{b}$ which is the closure of the 
second strand of $b$ (resp. which is the closure of the rest of the strand of $b$). 
By using the argument in  \cite[Section~11]{McMullen}, one sees that 
the Teichm\"{u}ler  polynomial $P(t_1, t_2, u)$ is given by 
$$P(t_1, t_2, u)= \det(uI- \sigma_2^{-1}(t_2) \sigma_1 (t_1) \sigma_2^{-1}(t_1)),$$ 
where 
$\sigma_1(t)= \left(\begin{array}{cc}t & t \\0 & 1\end{array}\right)$ and 
$\sigma_2^{-1}(t)= \left(\begin{array}{cc}1 & 0 \\t^{-1} & t^{-1}\end{array}\right)$. 
(Note  that our convention of the sign of braids is different from the one in \cite{McMullen}.)  
Hence 
$$P(t_1, t_2,u)= \tfrac{1}{t_2} - u - t_1 u - \tfrac{u}{t_2} - \tfrac{u}{t_1t_2}+ u^2.$$
Now we transform this to the polynomial using our basis. 
Let $\{\alpha^*, \beta^*, \gamma^*\}$ be the dual basis for $H_1(M_{\mathrm{magic}}; {\Bbb Z})$. 
We set $s_1= \alpha^*$, $s_2= \beta^*$, $s_3=\gamma^*$. 
By the construction of the homeomorphism $H$, one observes that 
$t_1= s_3$, $t_2= s_1 s_2^{-1} s_3^{-1}$ and $u=s_2$. 
We obtain 
$$P(s_1, s_2, s_3)= -s_1 -s_2 + s_3 + s_1 s_2 -s_1 s_3 -s_2 s_3.$$

\begin{thm}
\label{thm_Fried-Oertel-Poly}
The dilatation of $x \alpha + y \beta + z \gamma \in int(C_{\Delta}(\Bbb Z))$ is the largest real root of 
$$P(s_1= t^x, s_2=t^y, s_3=t^z)=  -t^x - t^y + t^z+ t^{x+y}-t^{x+z}-t^{y+z}.$$
In particular, the dilatation of $x \alpha+ y \beta \in int(C_{\Delta}(\Bbb Z))$ is the largest real root of 
$t^{x+y}-2(t^x+t^y)+1.$
\end{thm}

\noindent
{\it Proof.} 
 We identify $H_2(M_{\mathrm{magic}}, \partial M_{\mathrm{magic}})$ with $H^1(M_{\mathrm{magic}}, \partial M_{\mathrm{magic}})$. 
 By \cite[Section~1]{McMullen}, 
 the dilatation of $a= x \alpha + y \beta + z \gamma $  is equal to the largest real root of 
 \begin{eqnarray*}
 &&  P(s_1= t^{a(\alpha^*)}, s_2= t^{a(\beta^*)}, s_3=t^{a(\gamma^*)})= P(t^x, t^y, t^z)
 \\
 &=& -t^x - t^y + t^z+ t^{x+y}-t^{x+z}-t^{y+z}.
 \end{eqnarray*}
 This completes the proof. 
 $\Box$
 \medskip
 
 \noindent
 Since 
$$P(t^x, t^y, t^z)= t^z(t^{x+y-z}-t^x - t^y - t^{x-z}- t^{y-z}+1),$$
the dilatation  of $x\alpha+ y \beta + z \gamma$ is the largest real root of 
$$f_{(x,y,z)}(t)= t^{x+y-z}-t^x - t^y - t^{x-z}- t^{y-z}+1.$$

\begin{rem}
Oertel obtained similar polynomial with respect to each fiber face  as the Teichm\"{u}ler  polynomial, see \cite{Oertel}. 
\end{rem}

\begin{rem}
\label{rem_twin}
For $x \ge y >0$ and $z>0$, 
$$f_{(x,y,-z)}(t)= f_{(x+z, y+z, z)}(t)= t^{x+y+z}- t^{x+z} - t^{y+z}-t^x-t^y+1.$$
Thus, 
$$\lambda(x \alpha + y \beta - z \gamma)= \lambda((x + z) \alpha + (y+z) \beta + z \gamma).$$ 
\end{rem}

\begin{thm}
\label{thm_minimal-ray}
The homology class $\alpha+ \beta$ realizes the minimal normalized entropy with respect to $\Delta$, i.e, 
the ray through  $\alpha+ \beta$ attains the minimum of 
$\overline{\mathrm{ent}}(\cdot)=X_T(\cdot) \mathrm{ent}(\cdot):  int(C_{\Delta}) \rightarrow {\Bbb R}$. 
\end{thm}

\noindent
{\it Proof.} 
Clearly $x \alpha + y \beta + z \gamma$ is in $ int(C_{\Delta})$ if and only if $y \alpha + x \beta + z \gamma $ is in $  int(C_{\Delta})$. 
The equality  $\overline{\mathrm{ent}}(x \alpha + y \beta + z \gamma) = \overline{\mathrm{ent}}(y \alpha + x \beta + z \gamma)$ holds,  
since these classes  have the same entropy and the same Thurston norm. 
Thus if the minimum of $\overline{\mathrm{ent}}$ is realized by the ray though $x \alpha + y \beta + z \gamma$, then $x$ must equal $y$.  

On the other hand, 
$x \alpha + x \beta + \gamma $ and $(x-1) \alpha + (x-1) \beta - \gamma$ have  the same Thurston norm $2x-1$ if $x >1$. 
By Remark~\ref{rem_twin}, they have the same entropy. 
Thus,
$$\overline{\mathrm{ent}}(x \alpha + x \beta + \gamma) = \overline{\mathrm{ent}}((x-1) \alpha + (x-1) \beta - \gamma) .$$
Since $\overline{\mathrm{ent}}$ is constant on each ray,  
we have 
\begin{eqnarray*}
\overline{\mathrm{ent}}(x \alpha + x \beta + \gamma) &=& \overline{\mathrm{ent}}( \alpha +  \beta +x^{-1} \gamma), 
\\
\overline{\mathrm{ent}}((x-1) \alpha + (x-1) \beta - \gamma) &=& \overline{\mathrm{ent}}( \alpha + \beta -  (x-1)^{-1} \gamma). 
\end{eqnarray*} 
The minimal ray does not pass through both $ \alpha +  \beta +x^{-1} \gamma$ and $\alpha + \beta -  (x-1)^{-1} \gamma$, 
because the minimum is realized  by a unique ray. 
Since $x(>1)$  is arbitrary, the desired ray must pass through $\alpha+ \beta$.  
$\Box$

\subsection{Fiber surface of genus $0$}

Let  $x \alpha+ y \beta + z \gamma int(C_{\Delta})$ be an integral homology class. 
Recall that  $F_{(x,y,z)}$ is connected if and only if  $x \alpha+ y \beta + z \gamma$ is primitive. 
Since $\{\alpha, \beta, \gamma\}$ is a basis of $H_2(M_{\mathrm{magic}}, \partial M_{\mathrm{magic}}; {\Bbb Z})$, we see that 
$x \alpha+ y \beta + z \gamma$ is primitive if and only if  $\gcd(x,y,z)= 1$, where $\gcd(x,y,0)$ is defined to be $ \gcd(x,y)$. 
The topological type of $F_{(x,y,z)}$ can be determined by Lemma~\ref{lem_topological-type}. 
In this section, we find all the homology classes  whose minimal representatives are connected and of genus $0$. 

By (\ref{equation_int-delta}),  if $x \ge y$, then $x \alpha+ y \beta + z \gamma$ is in  $int(C_{\Delta}) $ 
if and only if $x \ge y > 0$ and $y >z$. 
In this section, we consider those  classes  for simplicity. 

\begin{lem}
\label{lem_genus0}
Let $x,y, z \in {\Bbb Z}$. 
Suppose that $x \ge y >0$, $ y>z$ and $\gcd(x,y,z)=1$. Then 
\begin{equation}
\label{equation_key-inequality}
\gcd(x,y+z)+ \gcd(y,z+x)+ \gcd(z,x+y)-x-y+z-2 \le 0.
\end{equation}
If $(x,y,z)$ satisfies  the equality in (\ref{equation_key-inequality}), then $z \ge 0$.  
\end{lem}

\noindent
{\it Proof.}  
Note that 
$$- \chi(F_{(x,y,z)}) = -(2-2g - \sharp (\partial F_{(x,y,z)})) = x+y-z,$$
where $g$ denotes the genus of $F_{(x,y,z)}$.  
Hence 
\begin{equation}
\label{equation_genus} 
g= \tfrac{x+y-z+2 -  \sharp (\partial F_{(x,y,z)}) }{2} \ge 0.
\end{equation}
By substituting $\sharp (\partial F_{(x,y,z)})= \gcd(x,y+z)+ \gcd(y,x+z)+ \gcd(z,x+y) $ for (\ref{equation_genus}), we have the desired inequality. 
Suppose that $x \ge y > 0> z= -z'$ $(z' >0)$. Then 
\begin{eqnarray*}
&\ &\gcd(x,y-z')+ \gcd(y,-z'+x)+ \gcd(-z',x+y)-x-y-z'-2 
\\
&\le& x+y+z'-x-y-z'-2
\\
&= &
-2.
\end{eqnarray*}
This completes the proof. 
$\Box$

\begin{prop}
\label{prop_genus0}
Let $x,y,z \in {\Bbb Z}$.  
Suppose that $x \ge y >z \ge 0$ and $\gcd(x,y,z)=1$. 
Then the  equality  of (\ref{equation_key-inequality}) holds for $(x,y,z)$ if and only if $(x,y,z)$ is either 
\begin{description}
\item[(1)] 
$z=0$ and $\gcd(x,y)=1$, 
\item[(2)] 
$(x,y,z)=(n+1,n,n-1)$  for $n \not\equiv 0 \pmod{3}$ and $n \ge 2$, or 
\item[(3)] 
$(x,y,z)= (2n+1,n+1,n)$  for $n \ge 1$. 
\end{description}
\end{prop}

\noindent
The following proof was shown to the authors by Shigeki Akiyama. 
\medskip

\noindent
{\it Proof.} 
The  equality  of (\ref{equation_key-inequality}) holds for $(x,y,0)$ if and only if $\gcd(x,y)=1$, and hence we may suppose that $x \ge y >z >0$ 
by Lemma~\ref{lem_genus0}. 
It is easy to see that if $(x,y,z)$ is of  either type (2) or type (3), then it satisfies the equality of (\ref{equation_key-inequality}). 
To prove the ``only if'' part, we first show: 
\begin{claim}
\label{claim_aki}
Let $x,y,z \in {\Bbb N}$. 
Suppose that  $\gcd(x,y,z)=1$. 
Then $$ \{\gcd(N,x), \gcd(N,y), \gcd(N,z)\}$$ 
is pairwise coprime, where $N= x+y+z$. 
\end{claim}

\noindent
{\it Proof of Claim~\ref{claim_aki}.} 
We set  $\gcd(\gcd(N,x), \gcd(N,y))=k$. 
Then $k$ is a divisor of three integers $N,x$ and $y$. 
It is also a divisor of $z (= N - x-y)$. 
Since $\gcd(x,y,z)=1$, the integer $k$ must be $1$. 
This completes the proof the claim. 
\medskip

\noindent
Notice that the inequality  of (\ref{equation_key-inequality}) is equivalent to the inequality 
\begin{equation}
\label{equation_key2} 
N - \gcd(N,x)- \gcd(N,y)- \gcd(N,z) \ge 2z -2. 
\end{equation}
For all $N \le 79$, one can check that the statement  of Proposition~\ref{prop_genus0} is valid. 
We may suppose that $N \ge 80$. 
Since $x > \frac{N}{3}$ and $z < \frac{N}{3}$, we have $x \ge \frac{N+1}{3}$ and $z \le \frac{N-1}{3}$. 
Hence  $x \ge 27$. 
Let $p,q$ and $r $ be natural numbers  so that $p>q>r$ and 
$$\{p,q,r\}= \{\gcd(N,x), \gcd(N,y), \gcd(N,z)\}.$$ 
By Claim~\ref{claim_aki}, $\{p,q,r\}$ is pairwise coprime, and $p,q$ and $r$ are divisors of $N$.  
Therefore $pqr \le N$. 
This shows that 
\begin{equation}
\label{equation_hyoka}
\tfrac{N-p-q-r}{N} \ge \tfrac{pqr-p-q-r}{pqr} = 1 -\tfrac{1}{qr}- \tfrac{1}{p}(\tfrac{1}{q}+ \tfrac{1}{r}).
\end{equation}

If $r \ge 2$, then $q \ge 3$ and $p \ge 5$. 
Hence 
\begin{eqnarray*}
N-p-q-r \ge  N(1-\tfrac{1}{qr}- \tfrac{1}{p}(\tfrac{1}{q}+ \tfrac{1}{r}))=  
\tfrac{5}{6}N (1-\tfrac{1}{p} ) \ge \tfrac{2}{3}N > 2z.
\end{eqnarray*}
Thus, no $(x,y,z)$  satisfies the equality of (\ref{equation_key2})  in this case. 

We may suppose that $r=1$. 
If $q \ge 4$, then 
$$\tfrac{N-p-q-r}{N} \ge \tfrac{3}{4}- \tfrac{5}{4p}$$ 
by (\ref{equation_hyoka}). 
Since $z \le \frac{N-1}{3}$, 
we obtain 
$$N - p-q-r -2z +2 \ge N ( \tfrac{3}{4}- \tfrac{5}{4p}) - \tfrac{2N-8}{3}= \tfrac{8}{3}+ \tfrac{N(p-15)}{12p}.$$ 
If $p >15$, then $\frac{8}{3}+ \frac{N(p-15)}{12p}>0$, which implies that 
no $(x,y,z)$ satisfies the equality of (\ref{equation_key2}) in this case. 
If $p \le 14$, then $q \le 13$. We have 
$N-2z= x+y-z > x>26= -1+14+13$. 
Thus 
$$N-p-q-1 \ge N-14-13-1 > 2z-2,$$
which implies that no $(x,y,z)$ satisfies the equality of (\ref{equation_key2}) in this case.

We may suppose that $q \le 3$. 
It is enough to consider the equality of (\ref{equation_key2})  in case $(q,r)= (3,1), (2,1)$ and $(1,1)$. 
Take $w \in \{x,y,z\}$ so that $p = \gcd(N,w)$. 

\begin{description}
\item[(1)]  Case $(q,r)= (2,1)$ or $(3,1)$. 

Then $N- \gcd(N,w)-2-1 = 2z -2$ or $N- \gcd(N,w)-3-1= 2z-2$. 
We set $\gcd(N,w)= \frac{N}{k}$ $(k>1)$. Then 
$N ( 1 - \frac{1}{k}) \le \frac{2(N-1)}{3}+2$. 
Since we assume that   $N \ge 80  >16$, $k$ must be $2$ or $3$. 
\begin{description}
\item[(i)] Case $k=2$. 

Then $x= \frac{N}{2}$. 
If $(q,r)= (2,1)$, then $(x,y,z)= ( \frac{N}{2}, \frac{N+2}{4}, \frac{N-2}{4})$. 
We set $n= \frac{N-2}{4}$. 
We obtain $(x,y,z)= ( 2n+1, n+1, n)$, and such $(x,y,z)$ satisfies the equality (\ref{equation_key2}). 
If $(q,r)= (3,1)$, then $(x,y,z)= (\frac{N}{2},   \frac{N+4}{4}, \frac{N-4}{4} )$. 
In this case, $\gcd (N, \frac{N \pm 4}{4}) = \gcd (4, \frac{N \pm 4}{4} )$, which is a divisor of $4$. 
This does not occur since  $(q,r)= (3,1)$. 

\item[(ii)] Case $k=3$. 

If $(q,r)= (2,1)$, then $z = \frac{N}{3}- \frac{1}{2}$, which can not be an integer. 
Let $(q,r)= (3,1)$. 
Then $z = \frac{N}{3}-1$. 
Since $\gcd(N,w)= \frac{N}{3}$, we see that $w= \frac{N}{3}$ or $\frac{2N}{3}$. 
If $w= \frac{2N}{3}$, then $(x,y,z)= ( \frac{2N}{3}, 1,  \frac{N}{3}-1)$. 
This is a contradiction since $y >z$. 
If $w= \frac{N}{3}$, then $(x,y,z)= ( \frac{N}{3}+1, \frac{N}{3},  \frac{N}{3}-1)$. 
We set $n= \frac{N}{3}$. 
Then $(x,y,z)= (n+1,n,n-1)$. 
If $n \equiv 0 \pmod{3}$, then $\gcd(3n,n \pm 1)=1$, which is a contradiction since $(q,r)= (3,1)$. 
Otherwise, such $(x,y,z)$ satisfies the equality of (\ref{equation_key2}). 
\end{description}
\item[(2)] Case $(q,r)= (1,1)$. 

Then $N- \gcd(N,w)-1-1= 2z-2$. We have $N- \gcd(N,w) = 2z < \frac{2N}{3}$. 
This implies that  $\gcd(N,w)= \frac{N}{2}$. 
Thus,  $(x,y,z)= ( \frac{N}{2}, \frac{N}{4}, \frac{N}{4})$, which does not occur since $y >z$. 
\end{description}
This completes the proof of Proposition~\ref{prop_genus0}. 
$\Box$
\medskip

\noindent
By Proposition~\ref{prop_genus0} and Lemma~\ref{lem_topological-type}, we immediately obtain the following 
which characterizes integral homology classes in $ int(C_{\Delta})$  whose minimal representatives are  
spheres with punctures. 

\begin{cor}
\label{cor_genus0} 
Let $x \alpha + y \beta + z \gamma \in  int(C_{\Delta})$ be an integral homology class. 
Suppose that $x \ge y$ and $\gcd(x,y,z)=1$. 
Then the genus of  $F_{(x,y,z)}$ is $0$ if and only if $(x,y,z)$ satisfies either  
\begin{description}
\item[(1)] 
$z=0$ and $\gcd(x,y)=1$, 
\item[(2)] 
$(x,y,z)=(n+1,n,n-1)$  for $n \not\equiv 0 \pmod{3}$ and $n \ge 2$, or 
\item[(3)] 
$(x,y,z)= (2n+1,n+1,n)$  for $n \ge 1$. 
\end{description}
In case (1), 
$$\sharp(\partial_{\alpha}F_{(x,y,z)}) = \sharp(\partial_{\beta}F_{(x,y,z)})=1\ \mbox{and}\ 
\sharp(\partial_{\gamma}F_{(x,y,z)}) =x+y.$$
In case (2), 
$$\{ \sharp(\partial_{\alpha}F_{(x,y,z)}), \sharp(\partial_{\gamma}F_{(x,y,z)})\}= \{1,3\} \ \mbox{and}\ 
\sharp(\partial_{\beta}F_{(x,y,z)})=n.$$
In case (3), 
$$\sharp(\partial_{\alpha}F_{(x,y,z)})=2n+1 \ \mbox{and}\ 
\{ \sharp(\partial_{\beta}F_{(x,y,z)}), \sharp(\partial_{\gamma}F_{(x,y,z)})\}= \{1,2\}.$$
\end{cor}

\noindent
Corollary~\ref{cor_genus0} implies that 
each mapping class $\phi_{(x,y,z)}= [\Phi_{(x,y,z)}]$ can be described by a braid  since $\Phi_{(x,y,z)}$ fixes one boundary component of $F_{(x,y,z)}$.

\begin{thm}
\noindent
\label{thm_asym-P}
\begin{description}
\item[(1)] $\displaystyle\lim_{n,m \to \infty} \overline{\mathrm{ent}}(n \alpha+ m \beta) = 2 \log(2+ \sqrt{3})$ if $\displaystyle\lim_{n,m \to \infty} \tfrac{n}{m}=1$. 
\item[(2)] $\displaystyle\lim_{n \to \infty} \overline{\mathrm{ent}}((n+1) \alpha + n \beta + (n-1) \gamma )= \infty$. 
\item[(3)] $\displaystyle\lim_{n \to \infty} \overline{\mathrm{ent}}((2n+1) \alpha+ (n+1) \beta + n \gamma) = \infty$. 
\end{description}
\end{thm}

\noindent
{\it Proof.} 
(1) We see that 
$\overline{\mathrm{ent}}(n \alpha+ m \beta) $ goes to $ \overline{\mathrm{ent}}(\alpha+ \beta)$ as 
$m,n$ go to $\infty$ with the condition $\displaystyle\lim_{n,m \to \infty} \tfrac{n}{m}=1$. 
By Theorem~\ref{thm_Fried-Oertel-Poly}, $\overline{\mathrm{ent}}(\alpha+ \beta)= 2\log(2+ \sqrt{3})$. 
\medskip
\\
(2) We see that 
$$\lim_{n \to \infty} \overline{\mathrm{ent}}(\tfrac{(n+1) \alpha+ n \beta + (n-1) \gamma}{n}) = \overline{\mathrm{ent}}(\alpha+ \beta + \gamma) =\infty$$ 
since $\alpha+ \beta + \gamma \in   \partial \Delta $. 
The proof of (3) is similar to the proof of (2). 
$\Box$

\subsection{Proposition~\ref{prop_minimal-homology}}

 Let $\mathcal{H}_{n}$ be the set of homology classes $x \alpha + y \beta + z \gamma \in int(C_{\Delta}(\Bbb Z))$, $x \ge y$ such that 
their  minimal representatives are $n$-punctured spheres. 
By Corollary~\ref{cor_genus0}, one can determine elements of $\mathcal{H}_n$. 
This section is devoted to prove: 

\begin{prop} 
\label{prop_minimal-homology}
The homology class which reaches the minimal dilatation among elements of $\mathcal{H}_{n}$ is as follows. 
\label{thm_main2}
\begin{description}
\item[(1)] $k \alpha+ (k-1) \beta$ in case $n= 2k+1$ for $k \ge 2$.
\item[(2)] 
$3 \alpha+ 2 \beta + \gamma$ in case $n=6$ and 
$(2k+1) \alpha + (2k-1) \beta$ in case  $n= 4k+2$ for $k \ge 2$.
\item[(3a)] 
$\alpha+ \beta$ in case $n=4$ and 
$(4k+3) \alpha + (4k-1) \beta$ in case  $n= 8k+4$ for $k \ge 1$.
\item[(3b)] 
$5 \alpha + 3 \beta + 2 \gamma$ in case $n=8$ and 
$(4k+5) \alpha + (4k+1) \beta$ in case $n= 8k+8$ for $k \ge 1$.
\end{description}
\end{prop}

\begin{lem}
\label{lem_prod-key}
Let $m>n>0$.
\begin{description}
\item[(1)] $\overline{\mathrm{ent}}((m+1) \alpha + m \beta + (m-1) \gamma)> \overline{\mathrm{ent}}((n+1) \alpha + n \beta + (n-1) \gamma)$. 
\item[(2)] $\overline{\mathrm{ent}}((2m+1) \alpha + (m+1) \beta + m \gamma)> \overline{\mathrm{ent}}((2n+1) \alpha + (n+1) \beta + n \gamma)$. 
\end{description}
\end{lem}

\noindent
{\it Proof.} 
Let us consider homology classes $x \alpha + y \beta + z \gamma$ with 
$(x,y,z)=  (\frac{n+1}{n+2}, \frac{n}{n+2}, \frac{n-1}{n+2})$ for each $n >0$. 
These classes are in the open face $ int(\Delta)$ and 
pass through the line $x = \frac{1}{3}t+1$, $y= \frac{2}{3}t+1$, $z= t+1$. 
Note that 
$$ \overline{\mathrm{ent}} ( (n+1) \alpha + n \beta +  (n-1) \gamma) 
= 
\mathrm{ent}((\tfrac{n+1}{n+2}) \alpha + (\tfrac{n}{n+2}) \beta +  (\tfrac{n-1}{n+2}) \gamma), $$
and it goes to $\infty$ as $n $ goes to $\infty$ by Proposition~\ref{thm_asym-P}(2). 
We have 
$$\tfrac{1}{3 \mathrm{ent}(2 \alpha + \beta + 0 \gamma)} \approx \tfrac{1}{2.887}> \tfrac{1}{4 \mathrm{ent}(3 \alpha + 2 \beta + \gamma)} \approx \tfrac{1}{2.931}.$$ 
Since  $\frac{1}{\mathrm{ent}(\cdot)}: int(C_{\Delta}({\Bbb Q})) \rightarrow {\Bbb R}$ is a strictly concave function, 
we have for all $m > n \ge 3$, 
$$ \tfrac{1}{4 \mathrm{ent}(3 \alpha + 2 \beta + \gamma)} 
 > \tfrac{1}{(n+2) \mathrm{ent}((n+1) \alpha + n \beta + (n-1) \gamma )}
 > 
\tfrac{1}{(m+2) \mathrm{ent}((m+1) \alpha + m \beta + (m-1) \gamma )}.
$$
This implies  (1). 
The proof of (2) is similar. 
$\Box$
\medskip

We set
\begin{itemize}
\item $\Delta_1= \{x \alpha + y \beta\ |\ x+y =1,\ x >0,\ y >0\} \subset int(\Delta)$ (see Figure~\ref{fig_face}(center)). 
\item $C_{\Delta_1} = \{x \alpha + y \beta\ |\  x >0,\ y >0\} $ (see Figure~\ref{fig_face}(right)).
\item  $C_{\Delta_1}(\Bbb {Z})  = \{a\ |\ a \in C_{\Delta_1}\ \mbox{is an integral class}\}$. 
\item $C_{\Delta_1}(\Bbb {Q}) = \{a\ |\ a \in C_{\Delta_1}\ \mbox{is a rational class}\}$. 
\end{itemize}

\begin{lem}
\label{lem_symmetry-entropy}
\noindent
\begin{description}
\item[(1)]  For $x,y \in {\Bbb N}$ such that $\gcd(x,y)=1$,  
the monodromy for $x \alpha + y \beta$ is conjugate to the inverse of the monodromy for $y \alpha+ x \beta$. 
\item[(2)] For $ x,y >0$,  we have $\mathrm{ent}(x \alpha + y \beta) = \mathrm{ent}(y \alpha+ x \beta)$. 
\end{description}
\end{lem}

\noindent
{\it Proof.} 
The existence of a $\pi$ rotation $f: (S^3, {\mathcal C}_3) \rightarrow (S^3, {\mathcal C}_3)$ with respect to the line $\mathtt{L}$ 
of Figure~\ref{fig_symmetry}(2) implies that the monodromy for $x \alpha + y \beta$ is conjugate to the one for $-y \alpha - x \beta$. 
This implies (1). 
The claim (2) is immediate from the expression for $f_{(x,y,z)}(t)$. 
$\Box$
\medskip

Fixing $n \in {\Bbb N}$, we set 
\begin{itemize}
\item 
$\Delta_n= \{x \alpha + y \beta \ |\ x>0,\ y>0,\ x+y= n\} \subset C_{\Delta_1}$. 
 

\item 
$\Delta_n(\Bbb N)= \{x \alpha+ y \beta \in \Delta_n\ |\ x,y \in {\Bbb N},\ \gcd(x,y)=1\}$. 
\end{itemize}

\begin{lem}
\label{lem_mini-real}
$\mathrm{ent}\big(\frac{n \alpha + n \beta}{2} \big) = \min\{\mathrm{ent}(a)\ |\ a \in \Delta_n\}$ for each $n \in {\Bbb N}$.
\end{lem}

\noindent
{\it Proof.} 
Recall that the restriction of $\frac{1}{\mathrm{ent}(\cdot)}$ on $\Delta_n$ is strictly concave. Note that 
$\mathrm{ent}(a) \to \infty$ as $a \to n \alpha$ or $n \beta$.  
Thus $\mathrm{ent}(\cdot)|_{\Delta_n} : \Delta_n \rightarrow {\Bbb R}$ has the unique minimum. 
By Lemma~\ref{lem_symmetry-entropy}(2), $\frac{n \alpha + n \beta}{2} $ attains the minimum. 
$\Box$
\medskip

\noindent

%

\begin{lem}
\label{lem_mini-monodromy}
For  $m \ge 3$, $\min\{\mathrm{ent}(a)\ |\ a \in \Delta_{m-1}(\Bbb N)\}$ is realized by: 
\begin{description}
\item[(1)] $ (k-1) \alpha + k \beta$ and $k \alpha + (k-1) \beta$ in case $m=2k$. 
\item[(2)] $(2k-1) \alpha+ (2k+1) \beta$ and  $(2k+1) \alpha + (2k-1) \beta$ in case $m= 4k+1 $. 
\item[(3a)] $\alpha + \beta$ in case $m=3$, and 
$(4k-1) \alpha + (4k+3) \beta$ and $(4k+3) \alpha + (4k-1) \beta$ in case $m=8k+3 $ for $k \ge 1$. 
\item[(3b)] $(4k+1) \alpha+ (4k+5 )\beta$ and $(4k+5) \alpha + (4k+1) \beta$ in case $m=8k+7$. 
\end{description}
 \end{lem}
 
\begin{figure}[htbp]
 \begin{center}
\includegraphics[width=2.8in]{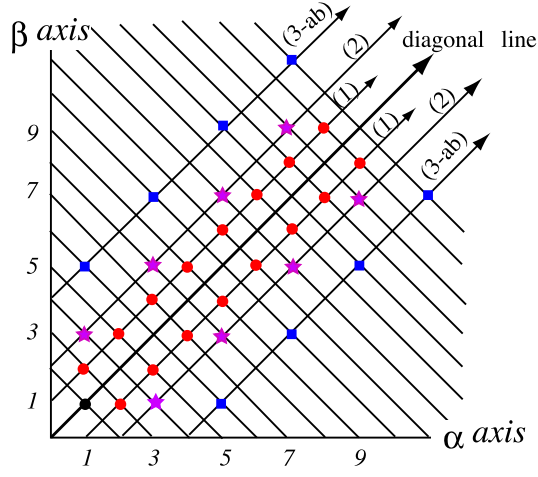}
\caption{the first quadrant of $\alpha \beta$ plane. homology classes in (1) (resp. (2), (3a, 3b) of Prop.~\ref{prop_mini-monodromy}) 
lie on the lines $\xrightarrow{\hspace{2mm}(1)}$  (resp. $\xrightarrow{\hspace{2mm}(2)}$ and $\xrightarrow{\hspace{2mm}(3-ab)}$).} 
\label{fig_line-plot}
\end{center}
\end{figure}

\noindent
{\it Proof.} 
The concavity of $\frac{1}{\mathrm{ent}(\cdot)}|_{\Delta_{n}}: \Delta_{n} \rightarrow {\Bbb R}$ and 
Lemma~\ref{lem_mini-real} tell us that 
if $|x-y|< |x'-y'|$ for $x \alpha + y \beta$, $x' \alpha+ y' \beta \in  \Delta_{n}$, then 
\begin{equation} 
\label{equation_useful}
\mathrm{ent}(x \alpha + y \beta) < \mathrm{ent}(x' \alpha+ y' \beta). 
\end{equation}
This implies that the minimal entropy among elements of $ \Delta_{m-1}(\Bbb N)$ is realized by 
$ (k-1) \alpha + k \beta$ or $k \alpha + (k-1) \beta$ if $m= 2k$. (See Figure~\ref{fig_line-plot}.)  
The proof for other cases can be shown in a similar way. 
$\Box$

\begin{lem}
\label{lem_monotonicity} 
For $k'> k >0$, we have the following. 
\begin{description}
\item[(1)]  $\overline{\mathrm{ent}}(k \alpha+ (k+c) \beta)> \overline{\mathrm{ent}}(k' \alpha+ (k'+c) \beta)$. 
\item[(2)]  $\mathrm{ent}(k \alpha+ (k+c) \beta)> \mathrm{ent}(k' \alpha+ (k'+c) \beta).$
\end{description}
\end{lem}

\noindent
{\it Proof.} 
For any $k>0$, we have 
$$\overline{\mathrm{ent}}(\tfrac{k \alpha + (k + c) \beta}{2k+c})=   
\mathrm{ent}(\tfrac{k \alpha + (k + c) \beta}{2k+c}) = (2k+c) \mathrm{ent}(k \alpha+ (k+c) \beta).$$
If $0 < k < k'$, then $\frac{c}{2k+c}> \frac{c}{2k'+c}$. 
By (\ref{equation_useful}), we see that 
$$\mathrm{ent}(\tfrac{k \alpha + (k + c) \beta}{2k+c}) > \mathrm{ent}(\tfrac{k' \alpha + (k' + c) \beta}{2k'+c}).$$ 
This implies  (1). 
By (1), 
\begin{eqnarray*}
&\ & \mathrm{ent}(\tfrac{k \alpha + (k + c) \beta}{2k+c}) = (2k+c) \mathrm{ent}(k \alpha+ (k+c) \beta)
\\
&>& \mathrm{ent}(\tfrac{k' \alpha + (k' + c) \beta}{2k'+c})=  (2k'+c) \mathrm{ent}(k' \alpha+ (k'+c) \beta)
\\
&>&  (2k+c) \mathrm{ent}(k' \alpha+ (k'+c) \beta). 
\end{eqnarray*}
Thus, $\mathrm{ent}(k \alpha+ (k+c) \beta)> \mathrm{ent}(k' \alpha+ (k'+c) \beta)$. 
This completes the proof of (2). 
$\Box$
\medskip

\noindent
{\it Proof of Proposition~\ref{prop_minimal-homology}.} 
(1) We consider the case $n=2k+1$.  
For $k=2$, we see that $\mathcal{H}_{5} = \{2 \alpha+ \beta\}$.  
If $k \ne 2$ and $2k \equiv 0 \pmod{3}$,  $\mathcal{H}_{2k+1}$ is the set of homology classes of type (1) of Corollary~\ref{cor_genus0}, that is 
$$\mathcal{H}_{2k+1}= \{x \alpha+ y \beta\ |\ x \alpha + y \beta \in \Delta_{2k-1}(\Bbb N), \ x \ge y\}.$$ 
In this case, $k \alpha + (k-1) \beta$ reaches the minimal entropy among elements of $\mathcal{H}_{2k+1}$ by Lemma~\ref{lem_mini-monodromy}(1). 
Otherwise (i.e,  $2 k \not\equiv 0 \pmod{3}$),  $\mathcal{H}_{2k+1}$ is the union of homology classes of type (1) and  (2) of Corollary~\ref{cor_genus0}: 
$$ \{(2k-2) \alpha + (2k-3) \beta + (2k-4)\gamma \} \cup \{x \alpha+ y \beta\ |\ x \alpha + y \beta \in \Delta_{2k-1}(\Bbb N), \ x \ge y\}.$$ 
One needs to compare the entropy for $(2k-2) \alpha + (2k-3) \beta + (2k-4)\gamma$ with the one for $k \alpha + (k-1) \beta$. 
In case $k =4$, $$\lambda(4 \alpha + 3 \beta) \approx 1.46557 < \lambda(6 \alpha+ 5 \beta + 4 \gamma) \approx 1.72208.$$
By Lemmas~\ref{lem_monotonicity} and \ref{lem_prod-key}, for $k  > 4$, we have 
\begin{eqnarray*}
&\ &(2k-1) \mathrm{ent}(k \alpha + (k-1) \beta)  
\\
 &\le& 7 \mathrm{ent}(4 \alpha +  3 \beta)
\\ 
&<& 7 \mathrm{ent}(6 \alpha+ 5 \beta + 4 \gamma)
\\
&<& (2k-1) \mathrm{ent}((2k-2) \alpha + (2k-3) \beta + (2k-4)\gamma).
\end{eqnarray*}
Thus, $ \mathrm{ent}(k \alpha + (k-1) \beta) < \mathrm{ent}((2k-2) \alpha + (2k-3) \beta + (2k-4)\gamma)$. 
This completes the proof.  
\medskip

\noindent
(2) Let us consider the case $n= 4k+2$.  
For $k=1$, 
$\mathcal{H}_{6} = \{3 \alpha+ \beta, 3 \alpha + 2 \beta + \gamma\}$. We have 
$$\lambda(3 \alpha + 2 \beta + \gamma) \approx 2.08102 < \lambda( 3 \alpha + \beta) \approx 2.29663.$$
For $k =2$,  $\mathcal{H}_{10} = \{7 \alpha + 4 \beta + 3 \gamma, 5 \alpha + 3 \beta, 7 \alpha+ \beta\}$. 
We have inequalities 
$$\lambda(5 \alpha + 3 \beta) \approx 1.41345< \lambda(7 \alpha + 4 \beta + 3 \gamma) \approx 1.55603\  \mbox{and}\ 
\lambda(5 \alpha + 3 \beta) < \lambda(7 \alpha+ \beta).$$ 
For $k = 3$, 
$\mathcal{H}_{14}= \{11 \alpha + 6 \beta + 5 \gamma, 7 \alpha + 5 \beta, 11 \alpha+ \beta, 11 \alpha + 10 \beta + 9 \gamma \}$. 
We have $\lambda(7 \alpha + 5 \beta)< \lambda(11 \alpha + \beta)$ by Lemma~\ref{lem_mini-monodromy}(2)  and  
\begin{eqnarray*}
\lambda(7 \alpha + 5 \beta) \approx 1.25141 &<& \lambda(11 \alpha + 6 \beta + 5 \gamma) \approx 1.39241 
\\
&<& \lambda(11 \alpha + 10 \beta + 9 \gamma) \approx 1.62913.
\end{eqnarray*}
By using the same arguments as in the case $n=2k+1$, we have 
for all $k >3$, $(4k)\mathrm{ent}((2k+1) \alpha + (2k-1) \beta)< 12 \mathrm{ent}(7 \alpha+ 5 \beta)$ and 
\begin{eqnarray*}
12 \mathrm{ent}(11 \alpha+ 6 \beta + 5 \gamma)&<& (4k) \mathrm{ent}((4k-1) \alpha+ 2k \beta + (2k-1)\gamma), 
\\ 
12 \mathrm{ent}(11 \alpha+ 10 \beta + 9 \gamma)&<& (4k) \mathrm{ent}((4k-1) \alpha+ (4k-2) \beta + (4k-3)\gamma). 
\end{eqnarray*} 
Thus, for all $k \ge 2$, $(2k+1) \alpha+ (2k-1) \beta$, which realizes the minimal entropy among elements of $\Delta_{4k}(\Bbb N)$, 
reaches the minimal entropy among elements of $\mathcal{H}_{4k+2}$. 
\medskip

 \noindent
 (3a) The proof for the  case $n= 8k+4$ is shown in a similar way. 
 \medskip

\noindent
(3b)  Let us consider the case $n= 8(k+1)$.  
For $k=0$,  we see that 
$\mathcal{H}_8= \{5 \alpha + 3 \beta + 2 \gamma, 5 \alpha + \beta, 5 \alpha + 4 \beta + 3 \gamma\}$, and 
$$\lambda(5 \alpha+ 3 \beta + 2 \gamma) \approx 1.72208 < \lambda(5 \alpha + 4 \beta + 3 \gamma) \approx 1.78164 < \lambda(5 \alpha+ \beta) 
\approx 2.08102.$$
For  $k \ge 1$, one can show that 
 $(4k+5) \alpha + (4k+1) \beta$ reaches the minimal dilatation among elements of $\mathcal{H}_{8k+4}$. 
$\Box$

\subsection{Monodromy}
\label{subsection_monodromies}

The braid $\Theta = \Theta_m= (\sigma_1^2  \sigma_2 \sigma_3 \cdots \sigma_{m-1})^{m-1} = 
(\sigma_1 \sigma_2 \cdots \sigma_{m-1})^m \in B_m$ 
is the full twist. Hence we have: 

\begin{figure}[htbp]
\begin{center}
\includegraphics[width=3.5in]{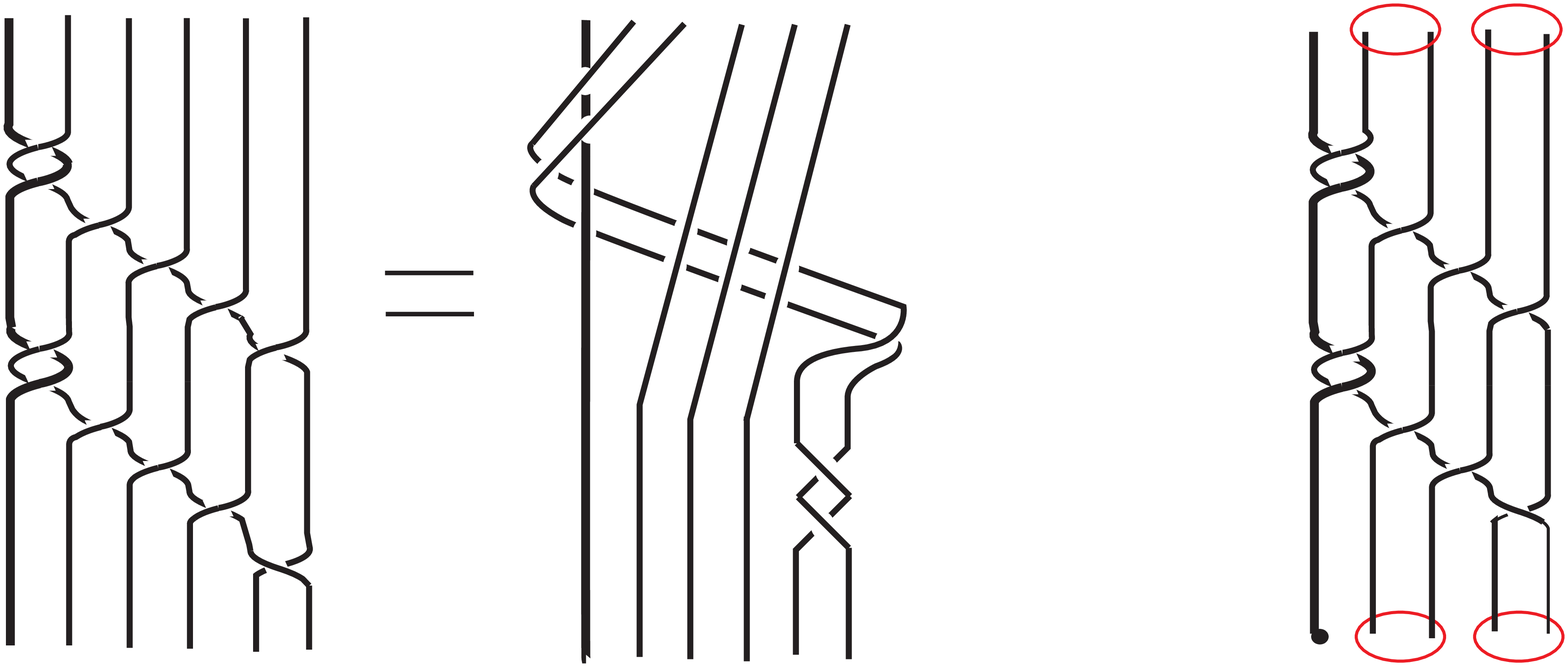}
\caption{(left) braid $T_{6,2}$.  (right) braid $T_{5,2}$.}
\label{fig_t_reducible}
\end{center}
\end{figure}

\begin{lem}
\label{lem_tmp-full}
If $p \equiv p' \pmod{m-1}$, then there exists an integer $k$ such that  
$T_{m,p} = T_{m,p'} \Theta^{k}$ and $\Gamma(T_{m,p})= \Gamma(T_{m,p'}) \in \mathcal{M}(D_m)$. 
\end{lem}


Let us consider the braid $T_{m,p}$ in case $\gcd(m-1,p) \ne 1$. 
For example, $T_{5,2}$ is a reducible braid, since a disjoint union of two simple closed curves in 
$D_5$ is invariant under $\Gamma(T_{5,2})$, see  Figure~\ref{fig_t_reducible}(right). 
It is not hard to see the following. 

\begin{lem}
\label{lem_tmp-reducible}
If  $\gcd(m-1,p) \ne 1$, then $T_{m,p}$ is a  reducible braid.  
\end{lem}

\begin{thm}
\label{thm_monodromy}
Suppose that $x \alpha + y \beta \in \Delta_{m-1}({\Bbb N})$ for $m \ge 3$. 
Then there exists $p= p(x,y)$ such that the braid $T_{m,p}$ is the 
monodromy on a fiber which is the minimal representative of $x \alpha + y \beta $. 
\end{thm}

\noindent
The rest of  section is devoted to proving Theorem~\ref{thm_monodromy} and explaining how to compute $p= p(x,y)$.

\subsubsection{Fiber surface}  
\label{subsubsection_Fiber-surface}

The aim of this section is to find fibers $F_{m(q),p(q)}$ for the magic manifold associated to sequences of natural numbers $q$ 
whose homology class $[F_{m(q),p(q)}] $ is in $ C_{\Delta_1} ({\Bbb Z})$. 

\begin{figure}[htbp]
\begin{center}
\includegraphics[width=5in]{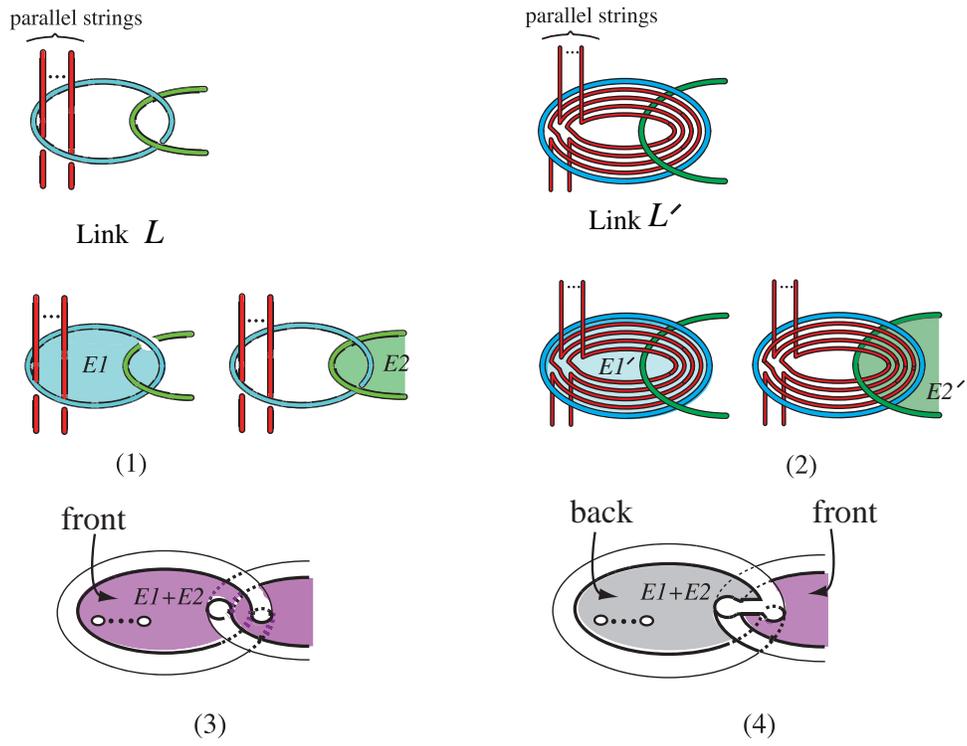}
\caption{(1) $E_1$, $E_2 \hookrightarrow E(L)$. (2) $E_1'$, $ E_2' \hookrightarrow E(L')$ (in case $N=2$ in Lemma~\ref{lem_xplusy}).  
(3) and (4)  $E_1+ E_2 \hookrightarrow E(L)$.}
\label{fig_xplusy}
\end{center}
\end{figure}

Let $L$ be a link in $S^3$. 
Let $E_1$ be an oriented disk with punctures which is   embedded  in the exterior $E(L) = S^3 \setminus \mathcal{N}(L)$ of $L$ and 
let $E_2$ be any embedded, oriented surface in $E(L)$ as in Figure~\ref{fig_xplusy}(1). 
The oriented surface $E_1+E_2$, which depends on  the orientation of $E_1$ and $E_2$,  is either of type (3) or type (4)  in  Figure~\ref{fig_xplusy}.  
The front  (resp. back) of $E_1+E_2$ is dark-colored (resp. light-colored) in the figure. 

Suppose that $E_1+E_2$ is of  type (3) (resp. (4)).  
Now, open $E(L)$ along $E_1$, and  let $E'$ and $E''$ be the resulting punctured disks obtained from $E_1$. 
Reglue $E'$ and $E''$ by 
twisting one of the disks by $360 \times N$ degrees in the clockwise (resp. counterclockwise) direction for some $N \in {\Bbb N}$. 
Then we obtain a new link, call it $L'$  such that  $E(L') \simeq E(L)$ 
(i.e, $E(L')$ is homeomorphic to  $E(L)$).  
Let $(E_1', E_2')$ be the ordered pair of  the embedded, oriented surfaces in $E(L)$ which are obtained from  the ordered pair $(E_1, E_2)$ 
(see Figure~\ref{fig_xplusy}(2)). 
The orientations of $E_1'$ and $E_2'$ are induced from $E_1$ and $E_2$ respectively.

\begin{lem}
\label{lem_xplusy} 
Let $L,L'$ be the links, $E_i, E'_i$ $(i =1,2)$ be the surfaces as above and let $N \in {\bf N}$ as above. 
There exists an orientation preserving homeomorphism $f: E(L) \rightarrow E(L')$ 
such that 
\begin{description}
\item[(1)] $f(yE_1+xE_2)= rE_1'+ xE_2'$ and
\item[(2)] $f(E_1) = E_1'$, 
\end{description}
where $y = xN + r$ for  $x,y \in {\Bbb N}$ and  $0 \le r <x$. 
In particular, 
\begin{description}
\item[(1')] $f(NE_1+ E_2)=  E_2'$. 
\end{description}
\end{lem}

\noindent
{\it Proof.} 
The construction  of $L'$ implies the existence of a homeomorphism $f: E(L) \rightarrow E(L')$ with the properties (1) and (2). 
(By using Figure~\ref{fig_xplusy}(3) and (4), one easily sees that  $f(E_1+E_2)= E_2'$ and $f(E_1) = E_1'$. 
One can generalize the first equality  to the claim (1).)  
$\Box$
\medskip

\noindent
Note that by Lemma~\ref{lem_xplusy}, 
$([E_1'], [E_2']) = ([E_1], [N E_1+ E_2])$. 

\begin{figure}[htbp]
\begin{center}
\includegraphics[width=3.8in]{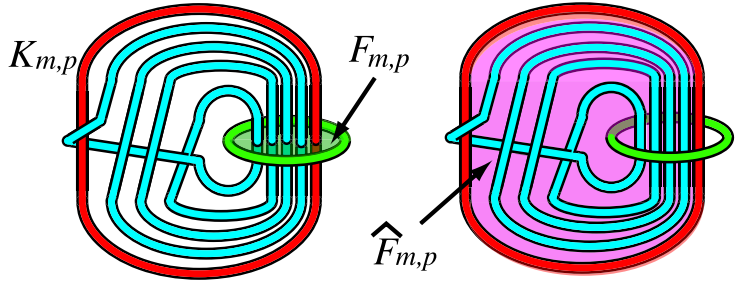}
\caption{$F_{m,p}$, $\widehat{F}_{m,p} \hookrightarrow E(\overline{T}_{m,p})$ in case $(m,p)= (5,1)$. 
(one sees that $\widehat{F}_{5,1} $ is a twice-punctured disk.)}
\label{fig_TwoSurfaces}
\end{center}
\end{figure}

Let us consider the exterior of the braided link $E( \overline{T}_{m,p})$. 
Now, we shall define two oriented surfaces  $ \widehat{F}_{m,p}$, $F_{m,p} \hookrightarrow E( \overline{T}_{m,p})$ whose orientations are induced by 
the oriented link $\overline{T}_{m,p}$. 
(Recall that  the orientation of $\overline{T}_{m,p}$ is given as in Figure~\ref{fig_artin}(right).) 
The oriented surface  $F_{m,p}$ is an $m$-punctured disk  which is bounded by the braid axis of $T_{m,p}$, see Figure~\ref{fig_TwoSurfaces}(left). 
Clearly, $F_{m,p}$ is a fiber for $E(\overline{T}_{m,p})$ with the monodromy $T_{m,p}$. 
The oriented surface  $ \widehat{F}_{m,p}$ is a $(p+1)$-punctured disk  which is bounded by  $K_{m,p}$, 
where $K_{m,p}$ is the knot which is  the closing the $1$st strand of $T_{m,p}$, see Figure~\ref{fig_TwoSurfaces}(right).

Given  $m,p \in {\Bbb N}$ and  $(k,\ell) \in {\Bbb N} \times {\Bbb N}$, 
the following construction  enables us to see another  fiber $F_{m', p'} \hookrightarrow E( \overline{T}_{m,p})$ 
with the monodromy $T_{m',p'}$. 
\medskip
\\
{\bf (Construction of fibers.)} 
\medskip
\\
Step 1. 
Apply Lemma~\ref{lem_xplusy} for the link $\overline{T}_{m,p}$, the ordered pair $(E_1, E_2)= (F_{m,p}, \widehat{F}_{m,p})$ and $N= k$. 
(Note that $F_{m,p}$ is a disk with punctures and hence one can apply Lemma~\ref{lem_xplusy}.)  
Then we obtain the ordered  pair of embedded surfaces in  $E( \overline{T}_{m, p})$; 
$$(F_{m, k(m-1)+p},  \widehat{F}_{m, k(m-1)+p})= (F_{m,p},  k F_{m,p} + \widehat{F}_{m,p}). $$
(Notice that  $\Theta^k= (\sigma_1^2 \sigma_2 \sigma_3 \cdots \sigma_{m-1})^{k(m-1)}$.) 
\medskip
\\
Step 2. 
Apply Lemma~\ref{lem_xplusy} for the link $ \overline{T}_{m, k(m-1)+p}$, the ordered pair 
$(E_1, E_2)=(\widehat{F}_{m, k(m-1)+p},  F_{m, k(m-1)+p})$ and $N= \ell$. 
(Note that $\widehat{F}_{m, k(m-1)+p}$ is a disk with punctures and hence one can apply the lemma.)  
Then it turns out that the new link $L'= L'( \overline{T}_{m, k(m-1)+p})$ 
is isotopic to  $\overline{T}_{m', p'}$, and we obtain 
 the ordered  pair of embedded surfaces   in  $E( \overline{T}_{m, p})$; 
\begin{equation}
 \label{equation_T2}
( \widehat{F}_{m',p'}, F_{m', p'}) 
=  (\widehat{F}_{m, k(m-1)+p}, \ell \widehat{F}_{m, k(m-1)+p}+F_{m, k(m-1)+p}),
\end{equation}  
where
\begin{eqnarray*}
m' &=&  \ell(k(m-1)+p)+m, 
\\
p'&=&  k(m-1)+p. 
\end{eqnarray*}
Thus we have another  fiber $F_{m', p'}$  for  $E(\overline{T}_{m,p}) $ with the monodromy $T_{m' p'}$.  
(End of the construction.) 
\medskip

\noindent
We sometimes denote $m'$ and $p'$  by $m'(m,p,k,\ell)$ and $p'(m,p,k,\ell)$. 
 
\begin{figure}[htbp]
 \begin{center}
\includegraphics[width=5in]{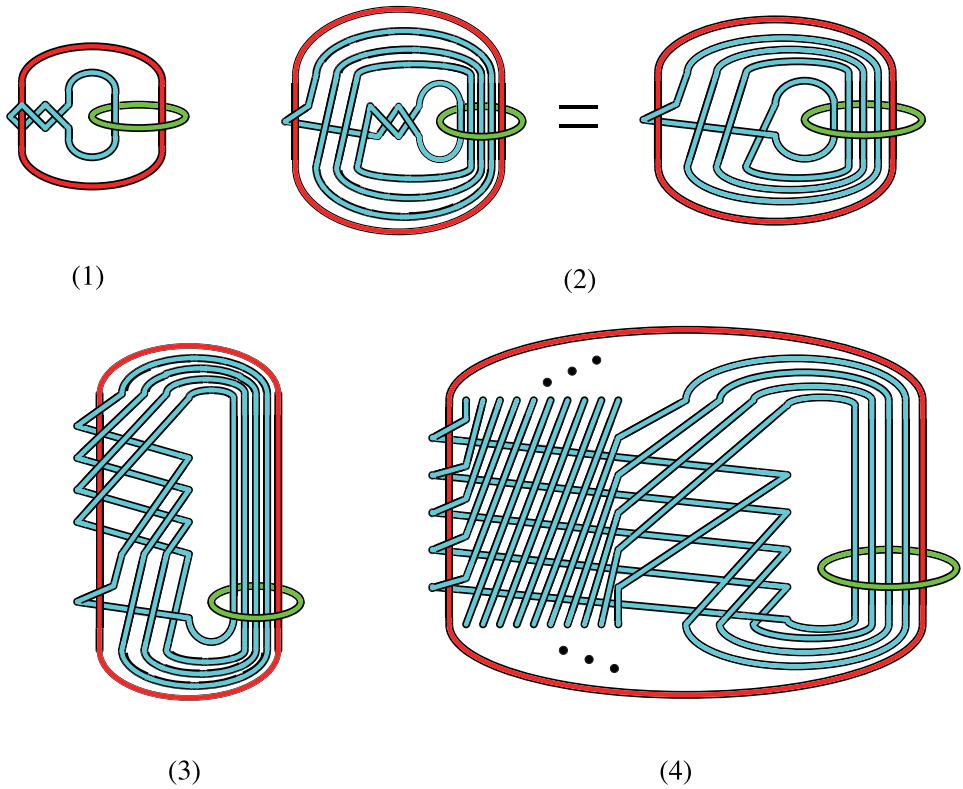}
\caption{construction of  fiber $F_{15,5}$ associated to $ (3,1,2)$. 
(1) $\mathcal{C}_3$. (2) $\overline{T}_{5,1}$. (fiber $F_{5,1}$ associated to $q= (3)$.) 
(3) $\overline{T}_{5,5}$. (4) $\overline{T}_{15,5}$.} 
\label{fig_T_15-5}
\end{center}
\end{figure}

 For example, in case $T_{m,p} = T_{5,1}$ and $(k,\ell) = (1,2)$, we have $T_{m', p'}= T_{15,5}$, see Figure~\ref{fig_T_15-5}(2),(3),(4). 
((3) explains Step 1 and (4) explains Step 2.) 

By using Lemma~\ref{lem_xplusy}, it is easy to see the following. 

\begin{prop}
\label{prop_T}
Let $x,y \in {\Bbb N}$. 
Suppose that $0 < x < y$  and $y \not\equiv 0 \pmod{x}$. 
Take $k,\ell, r_1, r_2$ such that 
\begin{eqnarray*}
y &=& xk + r_1\ (0 < r_1 < x,\ k \in {\Bbb N}), 
\\
x &=& r_1 \ell + r_2\ (0 \le  r_2 < r_1,\ \ell \in {\Bbb N}). 
\end{eqnarray*}
We apply the construction of a fiber  for a given $m,p \in {\Bbb N}$ and such a pair $(k,\ell)$. 
Then there exists an orientation preserving homeomorphism 
$f: E( \overline{T}_{m,p}) \rightarrow E( \overline{T}_{m', p'})$ 
such that 
$$f(x \widehat{F}_{m,p} + y F_{m,p}) = r_2  \widehat{F}_{m',p'}+ r_1 F_{m',p'} .$$
\end{prop}

Let $q= (q_1, q_2, \cdots, q_t)$ be a sequence of natural numbers. 
The number $t$ in the sequence, denoted by $|q| $, is called the {\it length} of $q$.  
For $q= (k_1, \ell_1, \cdots, k_j, \ell_j)$ with even length, 
let $q|_i=  (k_1, \ell_1, \cdots, k_i, \ell_i)$ for $i \le j$. 
For $q= (\ell_0, k_1, \ell_1, \cdots, k_j, \ell_j)$  with odd length, 
let $q|_i= (\ell_0, k_1, \ell_1, \cdots, k_i, \ell_i)$ for $i \le j$. 
Note that $q= q|_j$. 

\begin{figure}
\begin{center}
\includegraphics[width=3in]{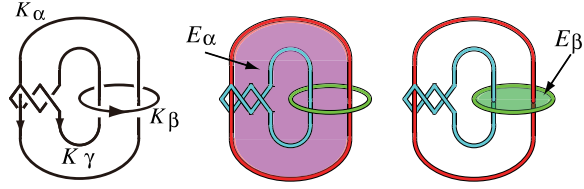}
\caption{(left) $\mathcal{C}_3$. 
(center), (right) $E_{\alpha}$, $E_{\beta} \hookrightarrow E( \mathcal{C}_3)$.}
\label{fig_Eab}
\end{center}
\end{figure}

\begin{figure}
\begin{center}
\includegraphics[width=5in]{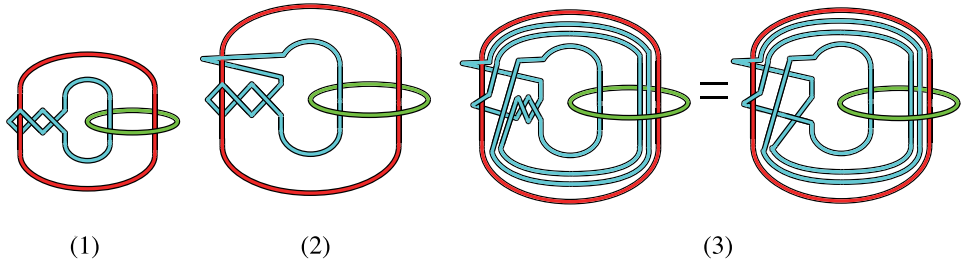}
\caption{fiber $F_{4,2}$ with monodromy $T_{4,2}$ associated to $(1,1)$. 
(2) and (3) describe the first and second part of Case 1 (even) respectively.} 
\label{fig_T_4-2}
\end{center}
\end{figure}

We will define  a fiber $F_{m(q),p(q)}$ for $E( {\mathcal C}_3)$ with the monodromy $T_{m(q),p(q)}$ associated to $q$ 
such that its homology class $[F_{m(q),p(q)}] $ is in $ C_{\Delta_1} ({\Bbb Z})$. 
To do so, we  define a fiber $F_{m(q|_i), p(q|_i)}$ for $E( {\mathcal C}_3)$ with the monodromy 
$T_{m(q|_i), p(q|_i)} $ inductively as follows. 
Another oriented diagram of  ${\mathcal C}_3$ is given in Figure~\ref{fig_Eab}(left). 
The oriented twice-punctured disk $E_{\alpha}$ (resp. $E_{\beta}$) bounded by $K_{\alpha}$ ($K_{\beta}$), 
whose orientation is induced by $\mathcal{C}_3$ 
is a  representative of  $\alpha$ (resp. $\beta$), see Figure~\ref{fig_Eab}(center, right). 
We first consider a sequence $q$ with even length. 
\medskip
\\
{\bf Case 1 (even).}   
Suppose that $q = (k_1, \ell_1)$.
First, apply Lemma~\ref{lem_xplusy} 
for $L = {\mathcal C}_3$, the ordered pair $(E_1, E_2)= (E_{\beta}, E_{\alpha})$ and $N= k_1$. 
Let $ (E_{\beta} , k_1 E_{\beta} + E_{\alpha})$ be the ordered pair of embedded surface in $E(\mathcal{C}_3) \simeq E(L')$
induced from $(E_{\beta}, E_{\alpha})$. 
Second, apply  Lemma~\ref{lem_xplusy} for $L'$, the ordered pair $ (k_1 E_{\beta} + E_{\alpha}, E_{\beta} )$ and $N= \ell_1$. 
Then we have the ordered pair of embedded surfaces 
\begin{equation}
\label{equation_c1e1order}
(k_1 E_{\beta} + E_{\alpha}, \ell_1(k_1 E_{\beta} + E_{\alpha}) + E_{\beta})
\end{equation}
in $E(\mathcal{C}_3) \simeq E(L'')$, 
where $L'' = (L')'$.  
We see that $L''$ is a braided link of $ T_{m(q), p(q)}= T_{(k_1+1)\ell_1+2, k_1+1}$, and 
\begin{equation}
\label{equation_c1e1}
(\widehat{F}_{m(q),p(q)}, F_{m(q),p(q)}) = (k_1 E_{\beta} + E_{\alpha}, \ell_1(k_1 E_{\beta} + E_{\alpha}) + E_{\beta})
\end{equation}
by (\ref{equation_c1e1order}). 
Therefore $ F_{m(q),p(q)}$ is a fiber for $E( {\mathcal C}_3)$ with the monodromy $T_{m(q), p(q)}$, and by (\ref{equation_c1e1}), 
$$\left[F_{m(q),p(q)}\right]= \ell_1 \alpha + (\ell_1k_1 +1) \beta \in C_{\Delta_1}(\Bbb Z)$$ 
since $\alpha= [E_{\alpha}]$ and $\beta = [E_{\beta}]$. 
For example in case $q= (1,1)$, we have 
$T_{m(q), p(q)}= T_{4,2}$, see Figure~\ref{fig_T_4-2}. 

Suppose that  $q= (k_1, \ell_1, \cdots, k_j, \ell_j)$, $j>1$. 
For $i=1$, we have defined a fiber $F_{m(q|_1),p(q|_1)}$ for $E( {\mathcal C}_3) $ with the monodromy 
$T_{m(q|_1), p(q|_1)}$ as above. 
Suppose that we have a fiber $F_{m(q|_i), p(q|_i)}$ for  $E( {\mathcal C}_3) \simeq E( \overline{T}_{m(q|_i), p(q|_i)})$  
with the monodromy $T_{m(q|_i), p(q|_i)}$. 
Apply the construction of a fiber  for $\overline{T}_{m(q|_i), p(q|_i)}$ and  the  pair $(k_{i+1}, \ell_{i+1})$. 
Then we have the ordered pair of embedded surfaces  $(\widehat{F}_{m(q|_{i+1}), p(q|_{i+1})},\  F_{m(q|_{i+1}),p(q|_{i+1})}) $ in $E( {\mathcal C}_3)$ 
(given in Step 2 in the construction) which is defined by 
\begin{eqnarray*}
(\widehat{F}_{m(q|_{i+1}), p(q|_{i+1})}, F_{m(q|_{i+1}),p(q|_{i+1})})
= (\widehat{F}_{m',p'}, F_{m',p'}), 
\end{eqnarray*} 
where 
$m'= m'(m(q|_i), p(q|_i), k_{i+1}, \ell_{i+1})$, 
$p'= p'(m(q|_i), p(q|_i), k_{i+1}, \ell_{i+1})$, see (\ref{equation_T2}). 
The surface 
$F_{m(q|_{i+1}), p(q|_{i+1})}$ is a fiber for $E( {\mathcal C}_3)$  with the monodromy $T_{m(q|_{i+1}), p(q|_{i+1})}$.  
By induction, it is shown that 
$\left[F_{m(q|_{i+1}), p(q|_{i+1})}\right] \in C_{\Delta_1}(\Bbb Z)$.

Next, let us consider a sequence $q$  with odd length.  
\medskip
\\
{\bf Case 2 (odd).} Suppose that $q= (\ell_0)$. 
Applying Lemma~\ref{lem_xplusy} for $ L= {\mathcal C}_3$, the ordered pair $(E_1,E_2)=(E_{\alpha}, E_{\beta})$ 
(\underline{not} $(E_1,E_2)=(E_{\beta}, E_{\alpha})$ as in Case 1 (even))  
and $N= \ell_0$, 
we obtain the ordered pair of embedded surfaces $(E_{\alpha}, \ell_0 E_{\alpha} + E_{\beta})$ in $E(\mathcal{C}_3) \simeq E(L')$.   
We see that $L'$ is a braided link of $T_{m(q),p(q)}=T_{\ell_0+2,1}$. 
We have 
$$(\widehat{F}_{m(q),p(q)} , F_{m(q),p(q)})= (E_{\alpha}, \ell_0 E_{\alpha} + E_{\beta}).$$ 
Therefore $ F_{m(q),p(q)}$ is a fiber for $E( {\mathcal C}_3)$ with the monodromy $T_{m(q), p(q)}$, and 
$$\left[ F_{m(q),p(q)}\right]= \ell_0 \alpha + \beta \in C_{\Delta_1}(\Bbb Z).$$
In case $\ell_0=3$, see Figure~\ref{fig_T_15-5}(2). 

Suppose that  $q= (\ell_0, k_1, \ell_1, \cdots, k_j, \ell_j)$, $2j+1>1$. 
For $i=0$, we have defined a fiber $F_{m(q|_0),p(q|_0)}$ for $E( {\mathcal C}_3) \simeq E(\overline{T}_{m(0), p(0)})$  
with the monodromy $T_{m(q|_0),p(q|_0)}$ as above. 
For $ i \ge 1$, in the same manner as Case 1 (even), 
a fiber $F_{m(q|_i), p(q|_i)}$ for $E( {\mathcal C}_3)$ with the monodromy $T_{m(q|_i), p(q|_i)} $ 
is given inductively, and we see that 
$\left[ F_{m(q),p(q)}\right] \in C_{\Delta_1}(\Bbb Z)$. 
In case $q=(3,1,2)$, see Figure~\ref{fig_T_15-5}.
\medskip

\subsubsection{Continued fraction}

Let us consider a continued fraction with length $j$ 
$$
w_1 + \frac{1}{w_2}{\genfrac{}{}{0pt}{}{}{+}}
           \frac{1}{w_3}{\genfrac{}{}{0pt}{}{}{+ \cdots +}}
           \frac{1}{w_{j-1}}{\genfrac{}{}{0pt}{}{}{+}}
           \frac{1}{w_{j}}{\genfrac{}{}{0pt}{}{}{}}:= 
w_1+ \cfrac{1}{w_2 + 
           \cfrac{1}{w_3+ \cdots
           \cfrac{1}{w_{j-1}+ 
           \cfrac{1}{w_{j}}}}}
$$
for  $w_i \in {\Bbb N}$.
We define $\left[w_1, w_2, \cdots, w_{j}\right] \in {\Bbb N}$ inductively as follows. 
\begin{eqnarray*}
\left[w_1\right]&=& w_1, 
\\
\left[ w_1 , w_2\right] &=& w_1 w_2 +1, 
\\
\left[w_1, w_2, \cdots, w_i \right] &=&   \left[w_1, w_2, \cdots, w_{i-1} \right]w_i+ \left[w_1, w_2, \cdots, w_{i-2}\right].    
\end{eqnarray*}
The following is elementary and well-known. 

\begin{lem}
\label{lem_elementary} 
\noindent
\begin{description}
\item[(1)]  $w_1 + \frac{1}{w_2}{\genfrac{}{}{0pt}{}{}{+}}
           \frac{1}{w_3}{\genfrac{}{}{0pt}{}{}{+ \cdots +}}
           \frac{1}{w_{j-1}}{\genfrac{}{}{0pt}{}{}{+}}
           \frac{1}{w_{j}}{\genfrac{}{}{0pt}{}{}{}} 
= \frac{\left[w_1, w_2, \cdots, w_{j}\right]}{\left[w_2, w_3, \cdots, w_{j}\right]}$.
\item[(2)] $\left[w_1, w_2, \cdots, w_{j}\right]= \left[w_{j}, w_{j-1}, \cdots, w_1\right]$. 
\end{description}
\end{lem}

\begin{definition}
\label{def_euclidean}
Suppose that $\gcd(u,v)= 1$ for $u,v \in {\Bbb N}$. 
We define two  sequences of non-negative integers 
$r = (  r_0, r_1, \cdots, r_{j+1})$ and $q = ( q_1, q_2,  \cdots, q_{j})$ associated to $\{u,v\}$ (according to the Euclidean algorithm). 
We set $r_0 = \max\{u,v\}$ and $r_1 = \min\{u,v\}$. 
Write  $r_{0} = r_{1} q_{1}+ r_2$ ($0 \le r_2 < r_{1}$). 
\begin{itemize}
\item 
If $r_{2} = 0$, then $r_{1}$ must be $1$ since $\gcd(u,v)=1$. 
We set 
$$r = (r_{0},r_{1}=1, r_2 =0)\ \mbox{and}\ q= (q_{1}).$$
\item 
Suppose that $r_2 \ne 0$. 
We define $q_2, q_3, \cdots$ and $r_3, r_4, \cdots$   inductively as follows. 
Let $q_i>0$ and $r_{i+1} \ge 0$ such that 
$r_{i-1}= r_i q_i + r_{i+1}$ ($0 \le r_{i+1} < r_i$).  
Since $r_0 >  \cdots > r_i > r_{i+1} \ge 0$, there exists  $j $ such that $r_{j+1}= 0$. 
(Then $r_{j}$ must be $1$ since $\gcd(u,v)=1$.) 
We set 
$$ r = ( r_0, r_1, \cdots, r_{j}=1, r_{j+1}=0)\ \mbox{and}\ q = ( q_0, \cdots, q_{j}).$$
\end{itemize}
\end{definition}
By using the sequence $q$, the fraction $ \frac{\max\{u,v\}}{\min\{u,v\}}$ 
can be expressed by the following two kinds of continued fractions. 
\begin{eqnarray}
\label{equation_continued-fraction}
\frac{\max\{u,v\}}{\min\{u,v\}}
&=&
q_1 + \frac{1}{q_2}{\genfrac{}{}{0pt}{}{}{+}}
           \frac{1}{q_3}{\genfrac{}{}{0pt}{}{}{+ \cdots +}}
           \frac{1}{q_{j-1}}{\genfrac{}{}{0pt}{}{}{+}}
           \frac{1}{q_{j}}{\genfrac{}{}{0pt}{}{}{}}, 
\\
\label{equation_continued-fraction2}
\frac{\max\{u,v\}}{\min\{u,v\}}
&=&
q_1 + \frac{1}{q_2}{\genfrac{}{}{0pt}{}{}{+}}
           \frac{1}{q_3}{\genfrac{}{}{0pt}{}{}{+ \cdots +}}
           \frac{1}{q_{j-1}}{\genfrac{}{}{0pt}{}{}{+}}
           \frac{1}{(q_{j}-1)}{\genfrac{}{}{0pt}{}{}{+}} 
             \frac{1}{1}{\genfrac{}{}{0pt}{}{}{}}     
\end{eqnarray}
with  length  $j$ and $j+1$ respectively. 
We  can choose the one with odd/even length among those continued fractions. 
\medskip

\noindent
{\it Proof of Theorem~\ref{thm_monodromy}.} 
Let $x \alpha + y \beta \in \Delta_{m-1}({\Bbb N})$. 
(By the definition of $ \Delta_{m-1}({\Bbb N})$, $x$ and $y$ are relatively prime.) 
From the continued fractions of $ \frac{\max\{x,y\}}{\min\{x,y\}}$ of the forms in (\ref{equation_continued-fraction}) and (\ref{equation_continued-fraction2}) 
(constructed by one of the sequence $(q_1, q_2, \cdots) = (w_1, w_2, \cdots)$ in Definition~\ref{def_euclidean} associated to $\{x,y\}$), 
we choose the one with odd length if  $x>y$ (resp. even length if $x<y$): 
\begin{equation}
\label{equation_suitable}
 \frac{\max\{x,y\}}{\min\{x,y\}}= 
w_1 + \frac{1}{w_2}{\genfrac{}{}{0pt}{}{}{+}}
           \frac{1}{w_3}{\genfrac{}{}{0pt}{}{}{+ \cdots +}}
           \frac{1}{w_{j-1}}{\genfrac{}{}{0pt}{}{}{+}}
           \frac{1}{w_{j}}{\genfrac{}{}{0pt}{}{}{}}. 
\end{equation}
Now, we take $s= (s_0, s_1, \cdots, s_{j+1})$ which is defined by 
\begin{eqnarray*}
s_0 &=& \max\{x,y\},
\\
s_1 &=& \min\{x,y\},
\\
s_{i+1} &=& s_{i-1} - s_i w_i \ \mbox{for\ } i \ge 1.
\end{eqnarray*}
Notice that $s_{j}= 1$ and $s_{j+1} = 0$. 
(If the continued fraction in (\ref{equation_suitable}) is of type (\ref{equation_continued-fraction}), then 
$s$ equals $r$ in Definition~\ref{def_euclidean}.) 

Suppose that $x<y$. 
(In this case, the continued fraction of (\ref{equation_suitable}) has even length.) 
Let us write 
$q= (w_1, w_2, \cdots, w_{j})= (k_1, \ell_1, \cdots, k_{j/2}, \ell_{j/2}) $. 
It is enough to show that 
a fiber $F_{m(q), p(q)}$ for $E(\mathcal{C}_3)$ associated to $q$ is a  representative of  $ x \alpha + y \beta $. 
(If this is the case, $F_{m(q), p(q)}$ is the minimal representative of $x \alpha + y \beta$ since $F_{m(q), p(q)}$ is a fiber.) 
By using Proposition~\ref{prop_T} repeatedly,  
we have 
\begin{eqnarray*}
 x \alpha + y \beta 
&=& \left[x E_{\alpha} + y E_{\beta}\right] (= \left[s_1 E_{\alpha} + s_0 E_{\beta}\right]) 
\\
&=& \left[s_3 \widehat{F}_{m(q|_1),p(q|_1)} + s_2 F_{m(q|_1),p(q|_1)} \right]
\\
&\vdots&  
\\
&=& \left[s_{j+1} \widehat{F}_{m(q|_{j/2}),p(q|_{j/2})} + s_{j} F_{m(q|_{j/2}),p(q|_{j/2})} \right] 
\\
&=&  \left[ 0 \widehat{F}_{m(q|_{j/2}),p(q|_{j/2})} + 1 F_{m(q|_{j/2}),p(q|_{j/2})} \right]  
\\
&=& \left[F_{m(q),p(q)}\right].
\end{eqnarray*} 
Since  the minimal representative of  $x \alpha + y \beta $ is an $(x+y+2)$-punctured sphere,  $m(q)$ equals $ x+y+1(=m)$. 

The proof for the case $x>y$ is similar. 
$\Box$

\subsubsection{Computation of $p=p(x,y)$ in Theorem~\ref{thm_monodromy}}

In this section we give a recipe to compute  $p$ in Theorem~\ref{thm_monodromy}. 
In Example~\ref{ex_observation-p}, we explain how the number $p$ is related to the pair $(x,y)$. 

Recall that $K_{m,p}$ is the knot obtained by the closing the $1$st strand of $T_{m,p}$. 
Let $K_{m,p}^{\diamondsuit}$ be the knot obtained by the closing the rest of strands, i.e, 
$K_{m,p}^{\diamondsuit}$ equals the closed braid of $T_{m,p}$ with $K_{m,p}$ removed. 
For  the braided link $\overline{T}_{m,p}$, we have a pair of  natural numbers  
$$(i(\widehat{F}_{m,p}, K_{m,p}^{\diamondsuit}), i(F_{m,p}, K_{m,p}^{\diamondsuit})) = (p,m-1),$$
where $i(S ,K )$ is the intersection number between the surface $S$ and the knot $K$. 
For $\mathcal{C}_3$, we have 
$$(i(E_{\alpha}, K_{\gamma}), i(E_{\beta}, K_{\gamma}))= (1,1),$$ 
see Figure~\ref{fig_Eab}.

\begin{ex}
\label{ex_observation-p}
 By the proof of Theorem~\ref{thm_monodromy} (see also 
 the argument in Case 2 (odd) in Section~\ref{subsubsection_Fiber-surface} and Figure~\ref{fig_T_15-5}), 
 $T_{15,5}$ is the monodromy on a fiber which is the minimal representative for $11 \alpha + 3 \beta$. 
 We explain  why   $p=5$ is derived from 
  $q= (q_1,q_2,q_3)=(3,1,2)$ and $r=(r_0,r_1,r_2,r_3, r_4)=(11,3,2,1,0)$ of Definition~\ref{def_euclidean} associated to $\{3,11\}$. 
We have 
$$\left[11 E_{\alpha}+ 3 E_{\beta}\right] = \left[2 \widehat{F}_{5,1}+ 3 F_{5,1}\right]= \left[2 \widehat{F}_{5,5}+ 1 F_{5,5}\right]= 
\left[0 \widehat{F}_{15,5}+ 1 F_{15,5}\right].$$ 
The following is a simple description of these equalities.  
\begin{equation}
\label{equation_simple-description}
\left.\begin{array}{ccccccc}
(r_0,r_1) &  & (r_2, r_1) &  & (r_2, r_3) &  & (r_4,r_3) \\
\parallel &  & \parallel &  & \parallel &  & \parallel \\
(11,3) & \xrightarrow[\ q_1=3\ ]{} & (2,3) &  \xrightarrow[\ q_2=1\ ]{} & (2,1) &  \xrightarrow[\ q_3=2\ ]{} & (0,1)  
 \end{array}\right.
\end{equation} 
In the process to find the fiber associated to $q=(3,1,2)$, we can find a sequence of pairs of intersection numbers 
$(1,1)$, $(1,4)$, $(5,4)$, $(5,14)$ 
obtained from 
$\mathcal{C}_3$, $\overline{T}_{5,1}$, $\overline{T}_{5,5}$, $\overline{T}_{15,5}$ respectively which is described from left to right as follows. 
\begin{equation}
\label{equation_intersection}
\left.\begin{array}{ccccccc}
(1,1) & \xleftarrow[\ q_1=3\ ]{} & (1,4) &  \xleftarrow[\ q_2=1\ ]{} & (5,4) &  \xleftarrow[\ q_3=2\ ]{} & (5,14)  
 \end{array}\right.
\end{equation}
Hence we can compute the number $p=5$ from the sequence $q= (3,1,2)$. 
To describe the number $p$ explicitly, 
we extend the sequence of (\ref{equation_simple-description}) to the left according to the Euclidean algorithm: 
\begin{eqnarray*}
\begin{array}{ccccccccc}
(r_0, r_0+r_1)& &(r_0,r_1) &  & (r_2, r_1) &  & (r_2, r_3) &  & (r_4,r_3) \\
\parallel &  &\parallel &  & \parallel &  & \parallel &  & \parallel \\
(11,14) &\xrightarrow[q_0=1]{} &(11,3) & \xrightarrow[q_1=3]{} & (2,3) 
&  \xrightarrow[q_2=1]{} & (2,1) &  \xrightarrow[q_3=2]{} & (0,1)  
\end{array}
\end{eqnarray*}
In the same way, 
we extend the sequence of  (\ref{equation_intersection}) to the left: 
\begin{eqnarray*}
\left.\begin{array}{ccccccccc}
(0,1)& \xleftarrow[\ q_0=1\ ]{} & (1,1) & \xleftarrow[\ q_1=3\ ]{} & (1,4) &  \xleftarrow[\ q_2=1\ ]{} & (5,4) &  \xleftarrow[\ q_3=2\ ]{} & (5,14)  
 \end{array}\right.
\end{eqnarray*}
These show that 
\begin{eqnarray*}
\frac{14}{11}&=&
1 + \frac{1}{3}{\genfrac{}{}{0pt}{}{}{+}}
         \frac{1}{1}{\genfrac{}{}{0pt}{}{}{+}}
           \frac{1}{2}{\genfrac{}{}{0pt}{}{}{}} 
= \frac{\left[1, 3,1,2\right]}{\left[3, 1,2\right]} = \frac{\left[q_0, q_1,q_2,q_3\right]}{\left[q_1, q_2,q_3\right]} , 
\\
\frac{14}{5}&=&
2 + \frac{1}{1}{\genfrac{}{}{0pt}{}{}{+}}
         \frac{1}{3}{\genfrac{}{}{0pt}{}{}{+}}
           \frac{1}{1}{\genfrac{}{}{0pt}{}{}{}} 
= \frac{\left[2, 1,3,1\right]}{\left[1, 3,1\right]} = \frac{\left[q_3, q_2,q_1,q_0\right]}{\left[q_2, q_1,q_0\right]} . 
\end{eqnarray*}
Thus the number $p(=5)$ in the question equals $\left[q_2, q_1,q_0\right]$. 
\end{ex}

\begin{prop}
\label{prop_computation-p}
Let $T_{m, p(x,y)}$ be the braid as in Theorem~\ref{thm_monodromy}. 
\begin{description}
\item[(1)]  
Let 
$ \frac{\max\{x,y\}}{\min\{x,y\}}= 
w_1 + \frac{1}{w_2}{\genfrac{}{}{0pt}{}{}{+}}
           \frac{1}{w_3}{\genfrac{}{}{0pt}{}{}{+ \cdots +}}
           \frac{1}{w_{j-1}}{\genfrac{}{}{0pt}{}{}{+}}
           \frac{1}{w_{j}}{\genfrac{}{}{0pt}{}{}{}}$
be the continued fraction chosen in (\ref{equation_suitable}). Then 
$$p= p(x,y)= [w_{j-1}, w_{j-2}, \cdots,w_1, w_0=1].$$
\item[(2)] $p= p(x,y)$ satisfies 
$$p \times \max\{x,y\} \equiv (-1)^{j} \pmod{x+y},$$
where $j$ is the length of the continued fraction of $ \frac{\max\{x,y\}}{\min\{x,y\}}$ in (1). 
\end{description}
\end{prop}

\noindent
{\it Proof.} 
(1) We have 
$$ \frac{x+y}{\max\{x,y\}}= w_0 + \frac{1}{w_1}{\genfrac{}{}{0pt}{}{}{+}}
           \frac{1}{w_2}{\genfrac{}{}{0pt}{}{}{+ \cdots +}}
           \frac{1}{w_{j-1}}{\genfrac{}{}{0pt}{}{}{+}}
           \frac{1}{w_{j}}{\genfrac{}{}{0pt}{}{}{}},$$
where $w_0=1$. 
It is not hard to show (1) by using the argument in Example~\ref{ex_observation-p}. 
\medskip
\\
(2)
By induction, one can show that 
\begin{eqnarray*}
&\ &
\left(\begin{array}{cc}w_0 & 1 \\1 & 0\end{array}\right) \left(\begin{array}{cc}w_1 & 1 \\1 & 0\end{array}\right) \cdots 
\left(\begin{array}{cc}w_j & 1 \\1 & 0\end{array}\right) \nonumber \\
&=& 
\left(\begin{array}{cc}\left[w_0, w_1, \cdots, w_{j}\right] & \left[w_0, w_1, \cdots, w_{j-1}\right] \\ 
\left[w_1, w_2, \cdots, w_{j}\right] & \left[w_1, w_2, \cdots,
      w_{j-1}\right] \end{array} \right). \ \ \ \ \ \ \ 
\end{eqnarray*}
Taking the determinant, one has 
\begin{eqnarray*}
(-1)^{j+1} &\equiv& - \left[w_0, w_1, \cdots, w_{j-1}\right] \left[w_1, w_2, \cdots, w_{j}\right] \pmod{\left[w_0, w_1, \cdots, w_{j}\right]}
\\
&\equiv& - \left[ w_{j-1}, w_{j-2}, \cdots, w_{0}\right] \left[w_1, w_2, \cdots, w_{j}\right]  \pmod{\left[w_0, w_1, \cdots, w_{j}\right]} .
\end{eqnarray*}
Note that 
$x+y = \left[w_0, w_1, \cdots, w_{j}\right] $, 
$p=   \left[w_{j-1}, w_{j-2}, \cdots, w_{0}\right]$, and 
$\max\{x,y\}= \left[w_1, w_2, \cdots, w_{j}\right] $. 
Thus, 
$$(-1)^{j+1} \equiv -p \times \max\{x,y\} \pmod{x+y}.$$
This implies (2).
$\Box$
\medskip

We show the converse of Theorem~\ref{thm_monodromy}. 

\begin{thm}
\label{thm_tmp-monodromy}
Suppose that $\gcd(p,m-1)=1$ for $p \ge 1$ and $m\ge 3$. 
Then there exist $x, y \in {\Bbb N}$ such that 
$T_{m,p}$ is the monodromy  on a fiber which is the minimal representative of $x \alpha + y \beta \in \Delta_{m-1}({\Bbb N})$. 
\end{thm}

\noindent
{\it Proof.}
Let $\varphi: {\Bbb N} \rightarrow {\Bbb N}$ be the Euler function. 
The number of braids $T_{m,p}$ satisfying $1 \le p \le m-1$ and $\gcd(p,m-1)=1$ equals $\varphi(m-1)$. 
Also, the number of elements $x \alpha + y \beta \in \Delta_{m-1}(\Bbb N)$  equals $\varphi(m-1)$. 
Let $x \alpha + y \beta$ and $ x' \alpha + y' \beta$ be distinct elements of $\Delta_{m-1}(\Bbb N)$. 
By Theorem~\ref{thm_monodromy}, 
it is enough to show that $p(x,y) \ne p(x', y') $ since we may assume that $1 \le p(x,y), p(x',y') \le m-1$ (see Lemma~\ref{lem_tmp-full}). 

Suppose that $(x,y) \ne (y',x')$. 
The concavity of $\mathrm{ent}(\cdot)|_{\Delta_{m-1}}: \Delta_{m-1} \rightarrow {\Bbb R}$ and Lemma~\ref{lem_symmetry-entropy} imply that 
$\mathrm{ent}(x \alpha + y \beta) \ne \mathrm{ent}(x' \alpha + y' \beta)$, and hence $T_{m,p(x,y)} \ne T_{m,p(x',y')}$ which implies that 
$p(x,y) \ne p(x',y')$. 

Suppose that $(x,y) = (y',x')$. 
(In this case, $ \mathrm{ent}(x \alpha + y \beta) = \mathrm{ent}(x' \alpha + y' \beta)$.) 
By Proposition~\ref{prop_computation-p}(2), we see that 
\begin{eqnarray}
\label{equation_p}
&\ &p(x,y) \times \max\{x,y\} + p(y,x) \times \max\{x,y\}  \nonumber \\ 
&=& (p(x,y)+ p(y,x)) \times \max\{x,y\}
\equiv 0 \pmod{x+y}.
\end{eqnarray}
Since $\gcd(\max\{x,y\}, x+y)=1$, we have $p(x,y)+ p(y,x) \equiv 0 \pmod{x+y}$. 
Thus,  $p(x,y) \not\equiv p(y,x) \pmod{x+y}$ which implies that 
$p(x,y) \ne p(y,x) (= p(x',y'))$. 
This completes the proof. 
$\Box$
\medskip

\noindent
Theorem~\ref{thm_tmp-monodromy} immediately gives: 

\begin{cor}
\label{cor_tmp-monodromy}
Suppose that $\gcd(p,m-1)=1$ for $p \ge 1$ and $m\ge 3$. 
Then  $ S^3 \setminus \overline{T}_{m,p}$  is homeomorphic to $S^3 \setminus {\mathcal C}_3$. 
\end{cor}

\begin{prop}
\label{prop_mini-monodromy}
Let  $m \ge 3$. 
The following shows homology classes realizing 
$\min\{\mathrm{ent}(a)\ |\ a \in \Delta_{m-1}(\Bbb N)\}$ and their monodromies. 
\begin{description}
\item[(1)] If $m= 2k$, then $ (k-1) \alpha + k \beta$ and $k \alpha + (k-1) \beta$ realize the minimum and their monodromies are given by 
$T_{2k,2}$ and $T_{2k,2k-3}$ respectively. 
\item[(2)] If $m= 4k+1 $, then  $(2k-1) \alpha+ (2k+1) \beta$ and  $(2k+1) \alpha + (2k-1) \beta$ realize the minimum and their monodromies are given by  
$T_{4k+1,2k+1}$ and  $ T_{4k+1,2k-1}$ respectively. 
\item[(3a)] If $m=3$, then $\alpha+ \beta$ realize the minimum and its monodromy is given by $T_{3,1}$. 
If $m=8k+3 $ $(k \ge 1)$, then $(4k-1) \alpha + (4k+3) \beta$ and $(4k+3) \alpha + (4k-1) \beta$ realize the minimum and their monodromies are given by 
$ T_{8k+3,2k+1}$ and  $T_{8k+3,6k+1}$ respectively.  
\item[(3b)] If $m=8k+7$, then $(4k+1) \alpha+ (4k+5 )\beta$ and $(4k+5) \alpha + (4k+1) \beta$ realize the minimum and their monodromies are given by 
$T_{8k+7,6k+5}$ and $T_{8k+7,2k+1}$ respectively.  
\end{description}
 \end{prop}
 
\noindent
{\it Proof.} 
We show the claim in case $m=2k$. 
Other cases can be shown in a similar way. 
By Lemma~\ref{lem_mini-monodromy}, 
the homology classes $ a = (k-1) \alpha + k \beta$ and $ a'=k \alpha + (k-1) \beta$ realize the minimum. 
Let us consider the monodromies $T_{m,p(k-1,k)}$ and $T_{m,p(k,k-1)}$. 
Let $(x,y)= (k-1,k)$. 
Since $x < y$, the continued  fraction which is chosen in (\ref{equation_suitable}) is 
$\frac{y}{x}= w_1+ \frac{1}{w_2}$, where $w_1=1$ and $w_2 = k-1$. 
By Proposition~\ref{prop_computation-p}(1), $p(k-1,k) = \left[ w_1, w_0 \right]=  \left[ 1, 1 \right]=2$. 
By (\ref{equation_p}), 
$$p(k-1,k)+ p(k,k-1) \equiv 0 \pmod{2k-1}.$$ 
Hence $p(k,k-1)= 2k-3$.  
By Lemma~\ref{lem_symmetry-entropy}(1), $T_{2k,2}$ or  $T_{2k,2k-3}$ gives the monodromy for $a$ and $a'$. 
$\Box$

\subsection{Proof of Theorem~\ref{thm_main}}

In Propositions~\ref{prop_minimal-homology} and \ref{prop_mini-monodromy}, we have proved Theorem~\ref{thm_main} except  $n= 6,8$. 
To complete the  proof, we shall  describe monodromies for two homology classes 
$3 \alpha + 2 \beta + \gamma$ and $5 \alpha + 3 \beta + 2 \gamma$ in Proposition~\ref{prop_two-monodromies}.

\begin{lem}
\label{lem_determination}
\noindent
\begin{description}
\item[(1)] The $5$-braided link $\overline{\sigma_1 \sigma_2^2 \sigma_3 \sigma_4}$ and the $4$-braided link 
$\overline{T}_{4,2}$ are isotopic to the $(-2,4,6)$-pretzel link. 
\item[(2)]  The braided link $\overline{b}$ for the $7$-braid $b$ as in Theorem~\ref{thm_main}(3b-i) 
 is isotopic to the $5$-braided link $\overline{\sigma_1 \sigma_2^2 \sigma_3 \sigma_4 \Theta_5^{-1}}$. 
\end{description}
\end{lem}

\begin{figure}[htbp]
 \begin{center}
\includegraphics[width=2.5in]{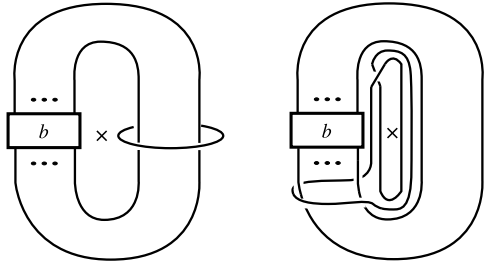}
\caption{(left) braided link $\overline{\beta}$. (right) closed braid representing  $\overline{\beta}$.} 
\label{fig_make-braid}
\end{center}
\end{figure}

\begin{figure}[htbp]
 \begin{center}
\includegraphics[width=5in]{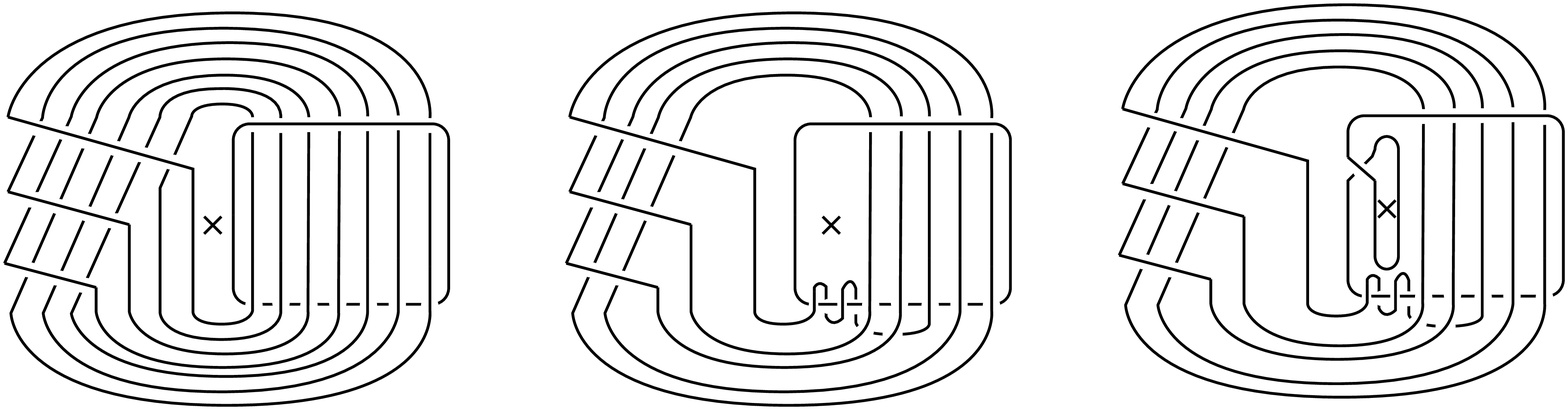}
\caption{(left)  braided link $\overline{b}$. (center) closed $6$-braid representing $\overline{b}$. (right) closed $7$-braid $\widehat{b'}$ 
representing $\overline{b}$.} 
\label{fig_8-exception}
\end{center}
\end{figure}

\noindent
{\it Proof.} 
(1) This is an easy exercise and we leave the proof for the readers. 
(Note: $T_{4,2}$ is conjugate to the $4$-braid $\sigma_1^2 \sigma_2 \sigma_3 \sigma_1^2 $, and it might be easier to see 
$\overline{\sigma_1^2 \sigma_2 \sigma_3 \sigma_1^2 }$ is isotopic to the  $(-2,4,6)$-pretzel link.)
\medskip
\\
(2) Let $\beta $ be an $n$-braid. By deforming the axis of $\beta$,  the braided link $\overline{\beta}$ can be represented 
by the closed braid $\widehat{\beta'}$ of $\beta' \in B_{n+2}$, 
where 
$\beta' = \sigma_{n+1}^{\varepsilon_1} \beta \sigma_n^{\varepsilon_2} \sigma_{n-1}^{\varepsilon_2} \cdots \sigma_1^{\varepsilon_2} \sigma_1^{\varepsilon_2} 
\sigma_2^{\varepsilon_2} \cdots \sigma_n^{\varepsilon_2}$ $(\varepsilon_1, \varepsilon_2 \in \{-1,1\})$, see Figure~\ref{fig_make-braid}. 
By using this method, 
$\overline{\sigma_1 \sigma_2^2 \sigma_3 \sigma_4 \Theta_5^{-1}}$ is represented by the closed $7$-braid $\widehat{a'}$, 
where 
$$a' = \sigma_6^{-1} (\sigma_1 \sigma_2^2 \sigma_3 \sigma_4 \Theta_5^{-1} ) \sigma_5^{-1} \sigma_4^{-1} \sigma_3^{-1} \sigma_2^{-1} \sigma_1^{-1} 
\sigma_1^{-1} \sigma_2^{-1} \sigma_3^{-1} \sigma_4^{-1} \sigma_5^{-1} .$$
On the other hand, the braided link $\overline{b}$ (Figure~\ref{fig_8-exception}(left)) can be represented by a closed $6$-braid 
as in Figure~\ref{fig_8-exception}(center) whose link type equals a closed $7$-braid as in Figure~\ref{fig_8-exception}(right). 
Namely,  $\overline{b}$ is isotopic to the closure of the $7$-braid $b'$: 
$$b'=  \underbar{6}\   \underbar{1}\   \underbar{2}\   \underbar{3}\   \underbar{4}\   \underbar{1}\   \underbar{2}\   \underbar{3}\   \underbar{1}\  
 \underbar{2} \  \underbar{5}^4\   \underbar{4} \  \underbar{3}\   \underbar{5}\   \underbar{4}\   \underbar{3}\   \underbar{2}\   \underbar{1}\   
  \underbar{1} \  \underbar{2}\   \underbar{3}\   \underbar{4}\   \underbar{5}, $$ 
  where \underbar{$i$} stands for $\sigma_i^{-1}$. 
We see that $a'$ is conjugate to  $b'$, 
since  the super summit set  for $a'$  is equal to the one for $b'$. 
(The super summit set is a complete conjugacy invariant, see \cite{EM}.) 
In fact, the super summit set consists of $4$ elements 
$\Theta_7^{-1}  1 2 3 4 3 2 1 5 4 3 6 5 4 3 2 1$, 
$\Theta_7^{-1}1 2 1 3 4 3 2 5 4 3 6 5 4 3 2 1$, 
$\Theta_7^{-1}  1 2 3 2 1 4 5 4 3 2 6 5 4 3 2 1 $ and 
$\Theta_7^{-1} 1 2 3 2 1 4 3 2 5 4 6 5 4 3 2 1 $, where 
$i$ stands for $\sigma_i$. 
(One can use the computer program ``Braiding'' by Gonz\'{a}lez-Meneses for a computation of the super summit set \cite{Gonzalez-Meneses}.)
Thus, the link types of  $\overline{b}$ and $\overline{\sigma_1 \sigma_2^2 \sigma_3 \sigma_4 \Theta_5^{-1}}$ are the same. 
This completes the proof. 
$\Box$

\noindent
Lemma~\ref{lem_determination} together with Corollary~\ref{cor_tmp-monodromy} implies: 

\begin{cor}\ 
\label{cor_determination}
\begin{description}
\item[(1)] 
$S^3 \setminus \overline{\sigma_1 \sigma_2^2 \sigma_3 \sigma_4}$ is homeomorphic to $S^3 \setminus \mathcal{C}_3$. 
\item[(2)] 
$S^3 \setminus \overline{b}$ is homeomorphic to $S^3 \setminus \mathcal{C}_3$. 
\end{description}
\end{cor}

\begin{prop}\ 
\label{prop_two-monodromies}
\begin{description}
\item[(1)] $\Gamma(\sigma_1 \sigma_2^2 \sigma_3 \sigma_4)$ is  the monodromy on a fiber 
which represents  $3 \alpha+ 2 \beta + \gamma$. 
\item[(2)] $\Gamma(b)$ is  the monodromy on a fiber which represents  $5 \alpha+ 3 \beta + 2\gamma$. 
\end{description}
\end{prop}

\noindent
{\it Proof.} 
(1) 
We see that  $\mathcal{H}_6= \{3 \alpha + 2 \beta + \gamma, 3 \alpha+ \beta\}$, see the proof of Proposition~\ref{prop_minimal-homology}. 
By Corollary~\ref{cor_genus0}, the monodromy for $3 \alpha + \beta$  permutes $4$ punctures cyclically and fixes two $1$ punctures. 
On the other hand, the monodromy for $3 \alpha+ 2 \beta + \gamma$ permutes $3$ punctures cyclically, 
and the mapping class $\Gamma(\sigma_1 \sigma_2^2 \sigma_3 \sigma_4)$ permutes $3$ punctures cyclically. 
By Corollary~\ref{cor_determination}(1), we complete the proof.  
\medskip
\\
(2)
We see that  $\mathcal{H}_8= \{5 \alpha + \beta, 5 \alpha + 3 \beta + 2 \gamma, 5 \alpha + 4 \beta + 3 \gamma\}$. 
The mapping class $\Gamma(b)$ permutes $5$ punctures cyclically, $2$ punctures cyclically  and fixes the other $1$ puncture. 
Among elements of $\mathcal{H}_8$, $5 \alpha + 3 \beta + 2 \gamma$ is the only class whose monodromy permutes $2$ punctures cyclically. 
By Corollary~\ref{cor_determination}(2),  we complete the proof. 
$\Box$

\section{Further discussion}
\label{section_further}

 \subsection{Pseudo-Anosov braids  with small dilatation}
\label{subsection_pA}

We consider the braids $T_{m,p}'$ defined in the introduction. 
The braid $T_{m,p}'$ may not be pseudo-Anosov, even though $T_{m,p}$ is so if $\gcd(p,m-1)=1$ (Corollary~\ref{cor_tmp-monodromy}). 
The inequality $\lambda(T_{m,p}') \le \lambda(T_{m,p})$ holds in case $T_{m,p}'$ is pseudo-Anosov. 
The following, which is clear by the definition of pseudo-Anosovs,  says when  the equality holds. 

\begin{lem}
\label{lem_entropy-equality} 
Suppose that $\gcd(p,m-1)=1$. 
Let $\Phi_{m,p}$ be the pseudo-Anosov homeomorphism which represents  $\Gamma(T_{m,p}) \in \mathcal{M}(D_m)$. 
Corresponding to the $1$st strand of $T_{m,p}$, there exists a puncture, say $a_{m,p}$, which is fixed by $\Phi_{m,p}$. 
Suppose that  the invariant foliation associated to $\Phi_{m,p} $ has no $1$-pronged singularity at $a_{m,p}$. 
Then $T_{m,p}'$ is  pseudo-Anosov such that 
$$\lambda(T_{m,p}') = \lambda(T_{m,p}).$$
\end{lem}

The families of braids $\{T'_{m,p}\}$ and $\{T_{m,p}\}$ contain examples with minimal dilatation. 
The following braids realize the minimal dilatation. 
\begin{itemize}
\item $T_{4,1}'  = \sigma_1 \sigma_2^{-1} \in B_3$, see Matsuoka \cite{Matsuoka}. 

\item $T_{5,1}' = \sigma_1 \sigma_2 \sigma_3^{-1}\in B_4$, see Ko-Los-Song \cite{KLS} and Ham-Song  \cite{HS}. 

\item $T_{6,2}' \sim \sigma_1 \sigma_2 \sigma_3 \sigma_4 \sigma_1 \sigma_2 \in B_5$, see Ham-Song \cite{HS}. 

\item $T_{6,3} \sim (\sigma_2 \sigma_1 \sigma_2 \sigma_1 (\sigma_1 \sigma_2 \sigma_3 \sigma_4 \sigma_5)^2)^{-1} \Theta \in B_6$, see 
Lanneau-Thiffeault  \cite{LT1}. 
 
 \item $T_{8,2}' \sim \sigma_4^{-2} (\sigma_1 \sigma_2 \sigma_3 \sigma_4 \sigma_5 \sigma_6)^2 \in B_7$, see Lanneau-Thiffeault  \cite{LT1}. 
  
 \item $T_{9,5}' \sim \sigma_2^{-1} \sigma_1^{-1} (\sigma_1 \sigma_2 \sigma_3 \sigma_4 \sigma_5 \sigma_6 \sigma_7)^5 \in B_8$, see Lanneau-Thiffeault \cite{LT1}. 
\end{itemize}
Here $b \sim b'$ means that $b$ is conjugate to $b'$.

All the braids in Proposition~\ref{prop_mini-monodromy} 
have been studied from the view point of their dilatations. 
Hironaka-Kin studied a family of braids 
$$\sigma_{(k)}= \sigma_1 \sigma_2 \cdots \sigma_{2k-2} \sigma_1 \sigma_2 \cdots \sigma_{2k-4} \in B_{2k-1}\hspace{2mm}(k \ge 3) $$
 with odd strands \cite{HK}. 
 It is easy to see that $\sigma_{(k)} \sim T_{2k,2}' $ (cf. Proposition~\ref{prop_mini-monodromy}(1)). 
Each braid $\sigma_{(k)} \in B_{2k-1}$ has the smallest known dilatation. 
Venzke found a family of braids $\{\psi_n\}$  with small dilatation \cite{Venzke}. 
\begin{eqnarray*}
\psi_n&=& L_n^2 \sigma_1^{-1} \sigma_2^{-1} \hspace{1.5cm}\mbox{if}\ n= 2k-1\ (k \ge 3), 
\\
\psi_n &=& L_n^{2k+1} \sigma_1^{-1} \sigma_2^{-1}  \hspace{1cm}\mbox{if}\ n=4k\ (k \ge 2), 
\\
\psi_n &=& L_n^{2k+1} \sigma_1^{-1} \sigma_2^{-1}  \hspace{1cm}\mbox{if}\ n=8k+2\ (k \ge 1), 
\\
\psi_n &=& L_n^{6k+5} \sigma_1^{-1} \sigma_2^{-1}  \hspace{1cm}\mbox{if}\ n=8k+6\ (k \ge 1), 
\\
\psi_6&=& \sigma_5 \sigma_4 \sigma_3 \sigma_2 \sigma_1 \sigma_5 \sigma_4 \sigma_3 \sigma_5 \sigma_4,
\end{eqnarray*}
where $L_n= \sigma_{n-1} \sigma_{n-2} \cdots \sigma_1 \in B_n$. 
It is not hard to see that 
$\psi_{2k-1} \sim T_{2k,2}'$,  
$\psi_{4k} \sim T_{4k+1,2k+1}' $, 
$\psi_{8k+2} \sim T_{8k+3,2k+1}' $, 
$\psi_{8k+6} \sim  T_{8k+7,6k+5}' $, and 
$\psi_6 \sim T_{6,2}$ (cf. Proposition~\ref{prop_mini-monodromy}(2)(3a)(3b)). 
By using Lemma~\ref{lem_sing-data} and Proposition~\ref{prop_mini-monodromy} together with Lemma~\ref{lem_entropy-equality}, we verify that
\begin{eqnarray*}
\lambda(\psi_{2k-1})&=&\lambda(T'_{2k,2}) =  \lambda(T_{2k,2}) ,
\\
\lambda(\psi_{4k}) &=&\lambda(T'_{4k+1,2k+1})  =  \lambda(T_{4k+1,2k+1}) ,
\\
\lambda(\psi_{8k+2}) &=&  \lambda(T'_{8k+3,2k+1})  = \lambda(T_{8k+3,2k+1}) , 
\\
\lambda(\psi_{8k+6}) &=&\lambda( T'_{8k+7,6k+5}) =  \lambda( T_{8k+7,6k+5}).
\end{eqnarray*}


Let $T_{(m)} \in B_m$ be either of  the two braids realizing the minimum in Proposition~\ref{prop_mini-monodromy}. 
For example, $T_{(2k)}= T_{2k,2}$ or $T_{2k,2k-3}$.  
Let $T_{(m)}' \in B_{m-1}$ be the braid obtained from $T_{(m)}$ by forgetting the $1$st strand of $T_{(m)}$. 
By using  Lemmas~\ref{lem_entropy-equality}, \ref{lem_sing-data} and Proposition~\ref{prop_mini-monodromy}, 
one has $\lambda(T_{(m)}) = \lambda(T_{(m)}')$. 
By Theorem~\ref{thm_Fried-Oertel-Poly} and Proposition~\ref{prop_mini-monodromy}, we have the following. 

\begin{cor}
\noindent
\label{cor_dilatation} 
\begin{description}
\item[(1)] 
$\lambda(T'_{(2k)}) $ equals  the largest real root of 
$$f_{(k-1,k,0)}(t)= t^{2k-1}- 2(t^{k-1}+ t^k)+1.$$
\item[(2)]  
$\lambda(T'_{(4k+1)})$ equals the largest real root of 
$$f_{(2k-1,2k+1,0)}(t)= t^{4k}-2(t^{2k-1}+ t^{2k+1}) +1.$$
\item[(3a)] 
$\lambda(T'_{(8k+3)})$ equals the largest real root of 
$$f_{(4k-1,4k+3,0)}(t)= t^{8k+2}-2(t^{4k-1}+ t^{4k+3})+1.$$
\item[(3b)] 
$\lambda(T'_{(8k+7)})$ equals the largest real root of 
$$f_{(4k+1,4k+5,0)}(t)= t^{8k+6}-2(t^{4k+1}+ t^{4k+5})+1.$$
\end{description}
\end{cor}

We now discuss the monotonicity of the dilatation of braids $T_{(m)}$. 
The following proposition is a corollary of Lemma~\ref{lem_monotonicity} and Proposition~\ref{prop_mini-monodromy}. 

\begin{prop}
\noindent
\label{prop_monotonicity-diagnal}
\begin{description}
\item[(1)] 
$\lambda(T_{(2k)})> \lambda(T_{(2(k+1))})$. 
\item[(2)]  
$\lambda(T_{(4k+1)})> \lambda(T_{(4(k+1)+1)})$. 
\item[(3a)] 
$\lambda(T_{(8k+3)})> \lambda(T_{(8(k+1)+3)})$. 
\item[(3b)] 
$\lambda(T_{(8k+7)})> \lambda(T_{(8(k+1)+7)})$.  
\end{description}
\end{prop}

\noindent
One can prove the following by using the argument in the proof of Lemma~\ref{lem_monotonicity}. 

\begin{lem}
\label{lem_monotonicity2} 
$\lambda(T_{(2k-1)}) > \lambda(T_{(2k)})$. 
\end{lem}

\noindent
In contrast to Lemma~\ref{lem_monotonicity2}, it is not true that 
$\lambda(T_{(2k)}) > \lambda(T_{(2k+1)})$ for all $k$. 
For example, 
\begin{eqnarray*}
\lambda(T_{(6)}) &<& \lambda(T_{(7)}), 
\\
\lambda(T_{(10)}) &<& \lambda(T_{(11)}). 
\end{eqnarray*}
See the computation of $\lambda(T_{(m)})$ and ${\mathrm{ent}}(T_{(m)})$ in the following table. 
We shall show $\lambda(T_{(2k)}) > \lambda(T_{(2k+1)})$ is true for other cases in the next.

\begin{tabular}{|c|c|c|c|}
\hline
$m$ & $T_{(m)}$  & $\lambda (T_{(m)})$ & ${\mathrm{ent}}(T_{(m)})$ \\ \hline
\hline
$3$ & $T_{3, 1}$  & $3.73205$ & $1.31696$ \\ \hline
$4$ & $T_{4, 2}$ or $T_{4, 1} $ & $2.61803$ & $0.962424$ \\ \hline
$5$ & $T_{5, 3}$ or $T_{5, 1} $ & $2.29663$ & $0.831443$ \\ \hline
$6$ & $T_{6, 2}$ or $T_{6, 3} $ & $1.72208$ & $0.543535$ \\ \hline
$7$ & $T_{7, 1}$ or $T_{7, 5} $ & $2.08102$ & $0.732858$ \\ \hline
$8$ & $T_{8, 2}$ or $T_{8, 5} $ & $1.46557$ & $0.382245$ \\ \hline
$9$ & $T_{9, 5}$ or $T_{9, 3} $ & $1.41345$ & $0.346031$ \\ \hline
$10$ & $T_{10, 2}$ or $T_{10, 7} $ & $1.34372$ & $0.295442$ \\ \hline
$11$ & $T_{11, 3}$ or $T_{11, 7} $ & $1.35293$ & $0.302271$ \\ \hline
$12$ & $T_{12, 2}$ or $T_{12, 9} $ & $1.27248$ & $0.240965$ \\ \hline
$13$ & $T_{13, 7}$ or $T_{13, 5} $ & $1.25141$ & $0.224273$ \\ \hline
$14$ & $T_{14, 2}$ or $T_{14, 11} $ & $1.22572$ & $0.203526$ \\ \hline
$15$ & $T_{15, 3}$ or $T_{15, 11} $ & $1.22257$ & $0.200958$ \\ \hline
$16$ & $T_{16, 2}$ or $T_{16, 13} $ & $1.19267$ & $0.176191$ \\ \hline
$17$ & $T_{17, 9}$ or $T_{17, 7} $ & $1.18129$ & $0.166609$ \\ \hline
$18$ & $T_{18, 2}$ or $T_{18, 15} $ & $1.16806$ & $0.155345$ \\ \hline
$19$ & $T_{19, 5}$ or $T_{19, 13} $ & $1.16432$ & $0.152136$ \\ \hline
$20$ & $T_{20, 2}$ or $T_{20, 17} $ & $1.14903$ & $0.13892$ \\ \hline
$21$ & $T_{21, 11}$ or $T_{21, 9} $ & $1.14192$ & $0.132708$ \\ \hline
$22$ & $T_{22, 2}$ or $T_{22, 19} $ & $1.13388$ & $0.125641$ \\ \hline
$23$ & $T_{23, 5}$ or $T_{23, 17} $ & $1.13071$ & $0.122845$ \\ \hline
$24$ & $T_{24, 2}$ or $T_{24, 21} $ & $1.12152$ & $0.114683$ \\ \hline
$25$ & $T_{25, 13}$ or $T_{25, 11} $ & $1.11665$ & $0.11033$ \\ \hline
$26$ & $T_{26, 2}$ or $T_{26, 23} $ & $1.11125$ & $0.105485$ \\ \hline
$27$ & $T_{27, 7}$ or $T_{27, 19} $ & $1.10869$ & $0.103176$ \\ \hline
$28$ & $T_{28, 2}$ or $T_{28, 25} $ & $1.10258$ & $0.0976543$ \\ \hline
$29$ & $T_{29, 15}$ or $T_{29, 13} $ & $1.09904$ & $0.0944354$ \\ \hline
$30$ & $T_{30, 2}$ or $T_{30, 27} $ & $1.09517$ & $0.0909069$ \\ \hline
$31$ & $T_{31, 7}$ or $T_{31, 23} $ & $1.09309$ & $0.0890074$ \\ \hline
$32$ & $T_{32, 2}$ or $T_{32, 29} $ & $1.08875$ & $0.0850323$ \\ \hline
$33$ & $T_{33, 17}$ or $T_{33, 15} $ & $1.08606$ & $0.0825554$ \\ \hline
$34$ & $T_{34, 2}$ or $T_{34, 31} $ & $1.08315$ & $0.0798714$ \\ \hline
$35$ & $T_{35, 9}$ or $T_{35, 25} $ & $1.08144$ & $0.0782958$ \\ \hline
$36$ & $T_{36, 2}$ or $T_{36, 33} $ & $1.07821$ & $0.0753015$ \\ \hline
$37$ & $T_{37, 19}$ or $T_{37, 17} $ & $1.07609$ & $0.0733366$ \\ \hline
$38$ & $T_{38, 2}$ or $T_{38, 35} $ & $1.07382$ & $0.0712265$ \\ \hline
$39$ & $T_{39, 9}$ or $T_{39, 29} $ & $1.07241$ & $0.0699047$ \\ \hline
\end{tabular}

\begin{lem}
\label{lem_monotone}
$\lambda(T_{(2k)}) > \lambda(T_{(2k+1)})$ for all $k \ge 2$ but $k= 3,5$. 
\end{lem}

\noindent 
The following is used for the proof of Lemma~\ref{lem_monotone}. 

\begin{lem}
\label{lem_Wsister-lem}
Let $x' > x >1$ and let $y'$ be the positive number such that 
$$X_T(x' \alpha + y' \beta)(= x'+y') = X_T(x \alpha + (x-1) \beta)+1 (=2x).$$ 
If $\lambda(x' \alpha + y' \beta)< \lambda(x \alpha+ (x-1) \beta)$, then 
$$\lambda((x'+ \tfrac{1}{2})\alpha + (y' + \tfrac{1}{2})\beta)< \lambda((x+\tfrac{1}{2}) \alpha+ (x-\tfrac{1}{2}) \beta).$$
\end{lem}

\noindent
{\it Proof.} 
One can show the claim by using the same argument as in \cite[Proposition~4.17]{KT1}. 
$\Box$ 
\medskip

\noindent
{\it Proof of Lemma~\ref{lem_monotone}.} 
One has $\lambda(3 \alpha + \beta) < \lambda(2 \alpha + \beta)$. 
This together with Lemma~\ref{lem_Wsister-lem} implies that 
$$
\lambda((2k+1) \alpha + (2k-1)\beta) < \lambda(2k \alpha + (2k-1) \beta) \hspace{2mm}  \mbox{for\ all}\ k \ge 1.
$$
One has another inequality 
$\lambda(9 \alpha + 5 \beta) < \lambda(7 \alpha + 6 \beta)$. 
Hence by Lemma~\ref{lem_Wsister-lem}, 
for all $k \ge 2$, one has 
\begin{eqnarray*}
\lambda((4k+3) \alpha + (4k-1) \beta)&<& \lambda((4k+1) \alpha + 4k \beta), 
\\
\lambda((4k+5) \alpha + (4k+1) \beta) &<& \lambda((4k+3) \alpha + (4k+2) \beta). 
\end{eqnarray*}
This together with Proposition~\ref{prop_mini-monodromy} completes the proof. 
$\Box$ 
\medskip

As a corollary of Lemmas~\ref{lem_monotonicity2} and ~\ref{lem_monotone} together with the equality 
$\lambda(T_{(m)}) = \lambda(T_{(m)}')$, one has: 

\begin{prop}\ 
\label{prop_compare}
\begin{description}
\item[(1)] 
$\lambda(T_{(2k-1)}) > \lambda(T'_{(2k)})$ for all $k \ge 2$.

\item[(2)]  
$\lambda(T_{(6)}) < \lambda(T_{(7)}')$ and $\lambda(T_{(10)}) < \lambda(T_{(11)}')$. 
For all $k \ge 2$ but $k= 3,5$, 
$\lambda(T_{(2k)}) > \lambda(T'_{(2k+1)})$. 
\end{description}
\end{prop}

\noindent
In particular, $T_{(10)} \in B_{10}$ has smaller dilatation than the Venzke's conjectural minimum $\lambda(\psi_{10}) (= \lambda(T_{(11)}'))$.

We turn to the asymptotic behavior of the normalized entropy of  the braid $T_{(m)}$. 
By Theorem~\ref{thm_asym-P}(1) and Proposition~\ref{prop_mini-monodromy}, we obtain the following. 

\begin{cor}
The normalized entropy of  $T_{(m)}$ goes to the minimal normalized entropy with respect to $\Delta$ as $m$ goes to $\infty$, i.e, 
$$\lim_{m \to \infty} \overline{\mathrm{ent}}(T_{(m)})= \overline{\mathrm{ent}}(\alpha + \beta)= 2 \log(2 + \sqrt{3}).$$
\end{cor}

Finally, we propose a conjecture on the minimal dilatation of braids of $\ell$ strands for $\ell \ge 9$. 

\begin{conj}\ 
\label{conj_minimal-entropy} 
\begin{description}
\item[(1)] 
The braid $T'_{(2k)}$ realizes the minimal dilatation among $(2k-1)$-braids for all $k \ge 5$. 

\item[(2)] 
The braid $T_{(10)}$ realizes the minimal dilatation among $10$-braids. 
The braid $T'_{(2k+1)}$ realizes the minimal dilatation among $2k$-braids for all $k \ge 6$.  
\end{description}
\end{conj}

\subsection{Asymptotic behavior of  entropy function} 
\label{subsection_Asymptotic}

We consider asymptotic behaviors of the entropy function for a family of homology classes in Proposition~\ref{prop_genus0}. 

\begin{thm}
\label{thm_asymp}
Let $x \alpha + y \beta \in C_{\Delta_1}$. 
\begin{description}
\item[(1)] $\displaystyle \lim_{x,y \to \infty} \mathrm{ent}(x \alpha + y \beta) =0$. 
\item[(2)] $\displaystyle\lim_{y \to \infty} \mathrm{ent}(x \alpha + y \beta) = \tfrac{\log2}{x}.$
\end{description}
\end{thm}

\noindent
Of course, $\displaystyle\lim_{x \to \infty} \mathrm{ent}(x \alpha + y \beta) = \tfrac{\log2}{y}$ by symmetry. 
\medskip

\noindent
{\it Proof.} 
(1) We may suppose that $x \le y$.  
By \cite[Theorem~3.5]{LO}, we have an inequality 
$$\mathrm{ent}(a+b) \le \min\{\mathrm{ent}(a), \mathrm{ent}(b)\}$$ 
for $a,b \in int(C_{\Delta}) $. 
Hence 
for all $\varepsilon>0$ so that $x - \varepsilon>0$ and for all $\delta>0$, 
$$\mathrm{ent}(x \alpha + (x+ \delta) \beta) \le \min\{ \mathrm{ent}((x-\varepsilon )\alpha + x \beta), \mathrm{ent}(\varepsilon \alpha + \delta \beta)\}.$$ 
Notice that $\mathrm{ent}(\varepsilon \alpha + \delta \beta)$ goes to $\infty$ as $\varepsilon$ goes to $0$. 
If one takes $\varepsilon>0$ sufficiently small, then one may assume that  
$$\mathrm{ent}(x \alpha + (x+ \delta) \beta) \le  \mathrm{ent}((x-\varepsilon) \alpha + x \beta).$$ 
Since $\mathrm{ent}(\cdot)$ is continuous, we have 
$\mathrm{ent}(x \alpha + (x+ \delta) \beta) \le  \mathrm{ent}(x \alpha + x \beta)$. 
Thus, 
$$ \lim_{x \to \infty}\mathrm{ent}(x \alpha + (x+ \delta) \beta) \le \lim_{x \to \infty}  \mathrm{ent}(x \alpha + x \beta) 
= \lim_{x \to \infty} \tfrac{1}{x}\mathrm{ent}(\alpha +\beta)= 0.$$ 
Since $\delta>0$ is arbitrary, the proof is completed. 
\medskip
\\
(2)  
By Theorem~\ref{thm_Fried-Oertel-Poly}, the dilatation of 
$x \alpha + y \beta + 0 \gamma \in  C_{\Delta_1}({\Bbb Z})$ is the largest real root of 
$$P(t^x, t^y,t^0) =P(t^x, t^y,1) = t^y R_x(t)+ (R_x)_*(t), $$ 
where $R_x(t) = t^x - 2$. 
By Lemma~\ref{lem_asymptotic-root1},  the largest real root of $P(t^x, t^y,1)$ converges to $2^{1/x}$, 
which is the unique real root of $R_x(t)$,  as $y \to \infty$.  
This claim  can be extended to homology classes of $C_{\Delta_1}({\Bbb Q})$, that is 
the dilatation of $x \alpha + y \beta \in C_{\Delta_1}({\Bbb Q})$ converges to $2^{1/x}$ as $y \to \infty$.  
Since the entropy function on 
$C_{\Delta_1}({\Bbb Q})$ can be extended to $C_{\Delta_1}$ uniquely, the proof is completed. 
$\Box$

\begin{prop}
\label{prop_golden-mean}
The entropy of $(n+1) \alpha + n \beta + (n-1)  \gamma \in int(C_{\Delta})$   
converges to the logarithm of the golden mean $\frac{1+ \sqrt{5}}{2}$ as $n$ goes to $\infty$. 
\end{prop}

\noindent
{\it Proof.} 
We have
$$P(t^{n+1}, t^n, t^{n-1}) = t^{n-1}\big(t^n(t^2-t-1)+ (t^2-t-1)_*\big).$$
If $(n+1) \alpha + n \beta + (n-1)  \gamma$ is an integral  class, then 
its dilatation $\lambda_n$ is the largest real root of $t^n(t^2-t-1)+ (t^2-t-1)_*$. 
The polynomial $t^2-t-1$ has the real root $\frac{1+ \sqrt{5}}{2}>1$. 
 By Lemma~\ref{lem_asymptotic-root1}, 
$\lambda_n$ converges to $\frac{1+ \sqrt{5}}{2}$ as $n \in {\Bbb N}$ goes to $\infty$. 
Since $\mathrm{ent}(\cdot)$ is continuous on $int(C_{\Delta})$, the proof is completed. 
$\Box$

\subsection{Relation between horseshoe braid and braid $T_{m,p}$} 

The horseshoe map was discovered by Smale around 1960. 
This map is well-known to be a simple factor possessing chaotic dynamics (\cite[Section~8.4.2]{Robinson} for example). 
For $\epsilon >0$, any $C^{1+\epsilon}$ surface diffeomorphism with positive topological entropy ``contains a horseshoe" in some iterate, see \cite{Katok} for more details.  
This tells us that the features of the horseshoe map is universal for  chaotic dynamical systems. 
In this section, we relate monodromies for homology classes in $C_{\Delta_1}({\Bbb Z})$ to the horseshoe map. 

\begin{figure}[htbp]
\begin{center}
\includegraphics[width=5in]{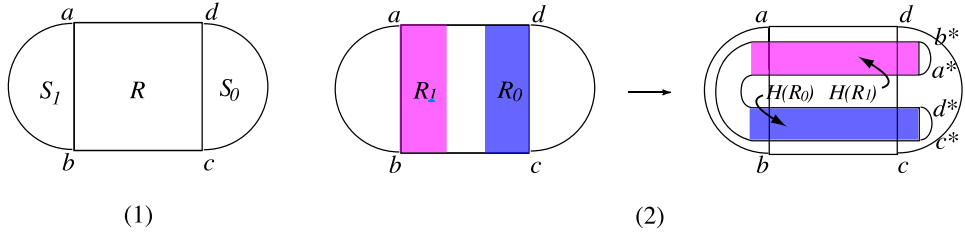}
\caption{(1) $R$, $S_0$, $S_1 \subset D$. (2) horseshoe map $H$. ($a^*$ is the image of $a$ under $H$, for example.)}  
\label{fig_S-horseshoe}
\end{center}
\end{figure}

The horseshoe map $H: D \rightarrow D$ is an orientation preserving diffeomorphism of the disk $D$ defined as follows. 
The action of $H$ on the rectangle $R$ and two half disks 
$S_0,S_1$ is given as in Figure~\ref{fig_S-horseshoe}. 
More precisely, 
the restriction $H|_{R_i}$ for $i \in \{0,1\}$ is an affine map such that 
$H$ contracts $R_i$ vertically and stretches horizontally, and $H|_{S_0 \cup S_1}: S_0 \cup S_1 \rightarrow S_0 \cup S_1$  is a contraction map. 
Then $H$ can be extended over the rest of $D$ without producing any new periodic points.

The set  $\Omega = \displaystyle\bigcap_{j \in {\Bbb Z}} H^j (R)$ is invariant under $H$. 
The map $H|_{\Omega}: \Omega \rightarrow \Omega$ can be described by using the symbolic dynamics as follows. 
We set $\mathcal{S}=\{0,1\}^{\Bbb Z} $, that is 
$\mathcal{S}$ is the the set of all two sided infinite sequences $s = (\cdots s_{-1} s_0 |s_1 \cdots)$ of $0$ and $1$, 
where we put the symbol $|$ between the $0$th element  and the $1$st element.  
We introduce the metric on $\mathcal{S}$ as follows. 
$$d(s,t) = \sum_{i \in {\Bbb Z}} \tfrac{|s_i - t_i|}{2^{|i|}}, $$
where $s =(\cdots s_{-1} s_0 |s_1 s_2 \cdots)$ and  $t = (\cdots t_{-1} t_0 |t_1 t_2  \cdots)$.

\begin{thm}[Smale]
\label{thm_horseshoe}
Let $\mathtt{s}: \mathcal{S} \rightarrow \mathcal{S}$ be the shift map, i.e, 
$\mathtt{s}$ is a homeomorphism such that 
\begin{eqnarray*}
\mathtt{s}(\cdots s_0| s_1 s_2 \cdots) = (\cdots s_0 s_1| s_2 \cdots).  
\end{eqnarray*}
The restriction  $H|_{\Omega}: \Omega \rightarrow \Omega$
is conjugate to the shift map   
$\mathtt{s}: \mathcal{S} \rightarrow \mathcal{S}$. 
The conjugacy $\mathcal{K}: \Omega \rightarrow \mathcal{S}$ is given by 
\begin{eqnarray*}
\mathcal{K}(x) &=& (\cdots \mathcal{K}_{-1}(x) \mathcal{K}_0(x) | \mathcal{K}_1(x) \cdots),
\hspace{2mm}
\mbox{where}
\end{eqnarray*}
\[
\mathcal{K}_j(x) =
\left\{
\begin{array}{ll}
1 \hspace{3mm}\  \mbox{if\ }&H^j(x) \in R_1,
\\
0 \hspace{3mm}\  \mbox{if\ }&H^j(x) \in R_0. 
\end{array}
\right.
\]
\end{thm}
If $x$ is a periodic point with the least period $k$ for $H$, then $\mathcal{K}(x)$ is a periodic sequence. 
The word $ \mathcal{K}_0(x) \mathcal{K}_1(x) \cdots \mathcal{K}_{k-1}(x)$ is called the {\it code} for $x$. 
Such word (modulo cyclic permutation) is said to be the {\it code} for the periodic orbit $\mathcal{O}_{H}(x) = \{x, H(x), \cdots, H^{k-1}(x)\}$.

\begin{rem}
\label{rem_horseshoe}
\noindent
\noindent
\begin{description}
\item[(1)]  Theorem~\ref{thm_horseshoe} asserts that there exists a one to one correspondence between the set of periodic points for $H|_{\Omega}$ and the 
set of periodic sequences in $\mathcal{S}$. 
\item[(2)] By using Theorem~\ref{thm_horseshoe}, one can show that the set of periodic points of $H|_{\Omega}$ is dense on $\Omega$. 
\end{description}
\end{rem}

\begin{figure}[htbp]
\begin{center}
\includegraphics[width=5in]{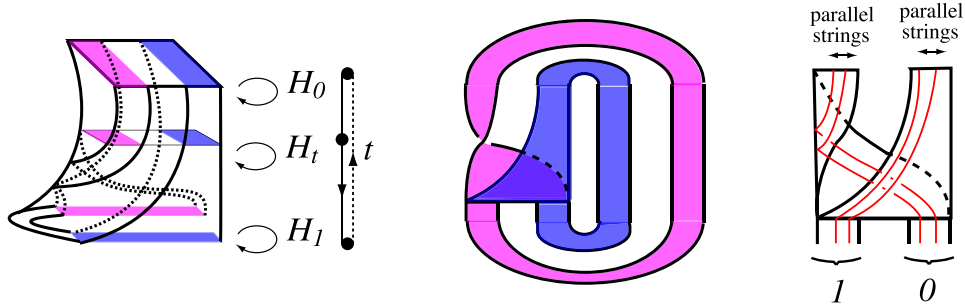}
\caption{(left) suspension of  horseshoe map. (center) horseshoe template $\mathcal{T}$. 
(right) non-horseshoe $4$-braid embedded in $\mathcal{T}$. (In this case, two strings of the braid are parallel.) }
\label{fig_template}
\end{center}
\end{figure}

Let $Q$ be a set of $n$ points consisting of periodic orbits of $H|_{\Omega}$.  
We take an isotopy $\{H_t\}_{t \in I = [0,1]}$ such that $H_0 = $ identity map on $D$ and $H_1 = H$. 
Then 
$$b(Q; \{H_t\}_{t \in I}) =\displaystyle\bigcup_{t \in I}H_t(Q) \times \{t\} \subset D \times I$$ 
is an $n$-braid. 
This depends on the choice of the isotopy, but it is determined uniquely up to a power of the full twist  $\Theta =(\sigma_1 \cdots \sigma_{n-1})^n$. 
Consider the suspension flow on the mapping torus by using a ``natural'' isotopy $\{H_t\}$, see Figure~\ref{fig_template}(left). 
For this isotopy, we denote the braid $b(Q; \{H_t\}_{t \in I}) $ by $b_{Q}$. 
By the definition of $H$, one can collapse  the image of the vertical lines of $R_0$ and $R_1$ under the isotopy to build the {\it horseshoe template} $\mathcal{T}$ as in Figure~\ref{fig_template}(center).  (For the template theory, see \cite{GHS}.) 
In this case 
the template is equipped with the semiflow induced by the suspension flow. 
It is easy to see that there exists a one to one correspondence between the set of periodic orbits of $H|_{\Omega}$ and the set of periodic orbits of the semiflow on $\mathcal{T}$. 
Each braid $b_Q$ can be embedded in $\mathcal{T}$ so that the closed braid of $b_{Q} $ becomes a finite union of periodic orbits of the semiflow on $\mathcal{T}$. 
Simply, we write $b_{Q}$ for the image of $b_{Q} \hookrightarrow \mathcal{T}$ when there exists no confusion.

Now, we define horseshoe mapping classes and horseshoe braids. 
Let $A_n$ be a set of $n$ points which lie on the horizontal line through the origin in the round disk $D$. 
We set an $n$-punctured disk $D_n= D \setminus A_n$. 
We say that $\phi \in \mathcal{M}(D_n)$ is a {\it horseshoe mapping class} if 
there exists a set of $n$ points $Q$ consisting of periodic orbits of $H|_{\Omega}$ and there exists an orientation preserving homeomorphism 
$g: D \setminus Q \rightarrow D_n$ such that $\phi$ is conjugate to  the mapping class $ [g \circ  H|_{D \setminus Q} \circ g^{-1}] \in  \mathcal{M}(D_n)$. 
A braid $\beta \in B_n$  is a {\it horseshoe braid} if the mapping class $\Gamma(\beta) \in \mathcal{M}(D_n)$  is a horseshoe mapping class. 
In other words, $\beta$ is a horseshoe braid 
if there exists an integer $k$ and  there exists a set of $n$ points 
consisting of a finite union of periodic orbits of $H|_{\Omega}$, denoted by $Q$, such that 
$\beta \Theta^k$ is conjugate to the braid $b_Q $. 
In this case,  there exists a braid $\gamma \in B_n$ such that 
$\gamma \beta \gamma^{-1} \Theta^k$ can be embedded in $\mathcal{T}$. 
However the converse is not true. 
For example, the $4$-braid of Figure~\ref{fig_template}(right) is not a horseshoe braid 
since there exists exactly one periodic orbit with the least period $2$ for $H|_{\Omega}$ whose code is $01$. 
By Remark~\ref{rem_horseshoe}(1), one can show that 
a braid $\beta$  embedded in $\mathcal{T}$ (ignoring the semiflow) is a horseshoe braid if  and only if no strings of the braid are  parallel. 
(See Figure~\ref{fig_template}(right).) 

\begin{prop}
\label{prop_horseshoe-type}
Suppose that $\gcd(p,m-1) =1$.  
If $1<p \le \frac{m-1}{2}$, then $T_{m,p} \in B_m$ is a horseshoe braid. 
\end{prop}

\begin{figure}[htbp]
\begin{center}
\includegraphics[width=5in]{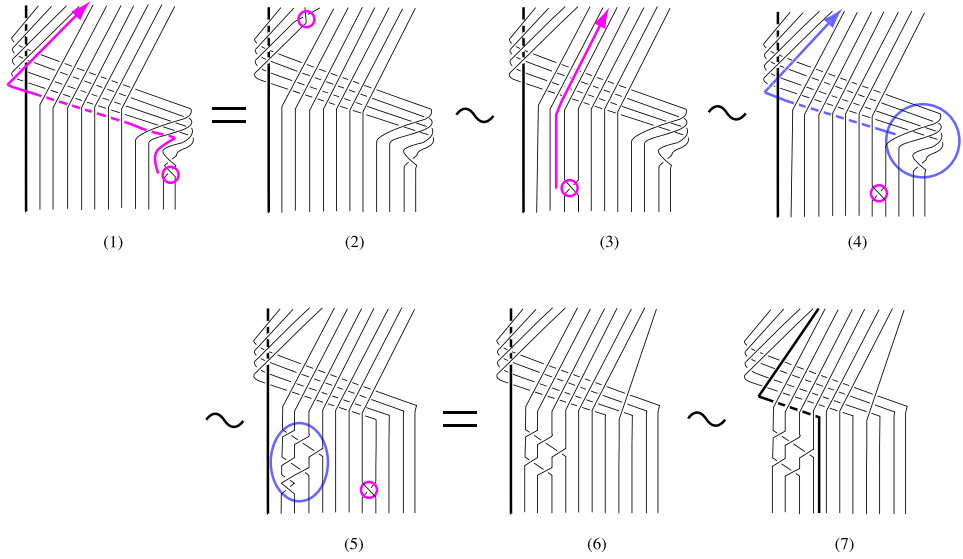}
\caption{conjugate braid of $T_{12,4}$. 
(1) $T_{12,4}$. (7) braid embedded in $\mathcal{T}$. }
\label{fig_T_12-4}
\end{center}
\end{figure}

\noindent
Obviously, if the braid $b$ is written by $b= c b'$, then $b$ is conjugate to $b'c$. 
This is used for the proof of Proposition~\ref{prop_horseshoe-type}. 
Before proving the proposition, we first see that $T_{12,4}$ is a horseshoe braid by using Figure~\ref{fig_T_12-4}. 

\begin{ex}
\label{ex_T_12-4}
The first braid of Figure~\ref{fig_T_12-4} is a representative of $T_{12,4}$. 
We slide the last crossing in the small circle to the top, see the second braid. 
Then it is conjugate to the third braid  of Figure~\ref{fig_T_12-4}. 
We repeat to slide  the last crossing in the small circle of the third braid  to the top. 
We see that 
it is conjugate to the fourth braid. 
The crossings in the large circle of the fourth braid can slide to the top, and then we see that the fourth braid is 
conjugate to the fifth braid which is isotopic to the sixth braid. 
Finally, it is easy to see that the sixth braid is conjugate to the seventh braid which can be embedded in $\mathcal{T}$. 
(In fact, the braid $\sigma_1 \sigma_2 \sigma_3 \sigma_4$ is a conjugacy.) 
Since no strings of the latter braid are parallel, one concludes that $T_{12,4}$  is a horseshoe braid. 
\end{ex}

\begin{figure}[htbp]
\begin{center}
\includegraphics[width=3.8in]{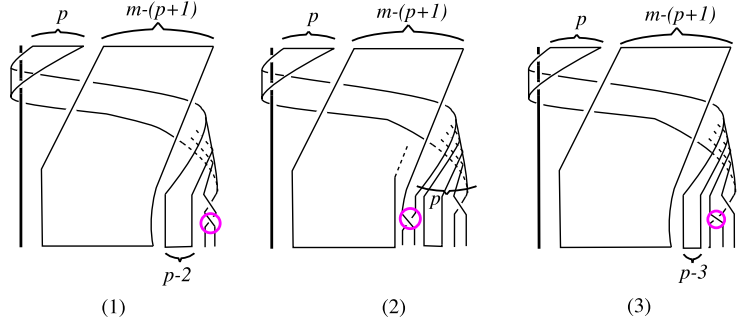}
\caption{(1) $T_{m,p}$. ($T_{m,p}$ is conjugate to either the braid drawn in (2) or the one drawn in (3).)}
\label{fig_horse-braid}
\end{center}
\end{figure}

\noindent
{\it Proof of Proposition~\ref{prop_horseshoe-type}.} 
We consider a representative of $T_{m,p}$ as in the first braid of Figure~\ref{fig_horse-braid}. (See also Figure~\ref{fig_t_reducible}(left).)  
By using the slide technique  in Example~\ref{ex_T_12-4}, we see that 
$T_{m,p}$ is conjugate to the second braid or the third  of Figure~\ref{fig_horse-braid}.  
(For example, $T_{12,4}$ is conjugate to the second type and $T_{12,5}$ is conjugate to the third type.) 

\begin{figure}[htbp]
\begin{center}
\includegraphics[width=4.5in]{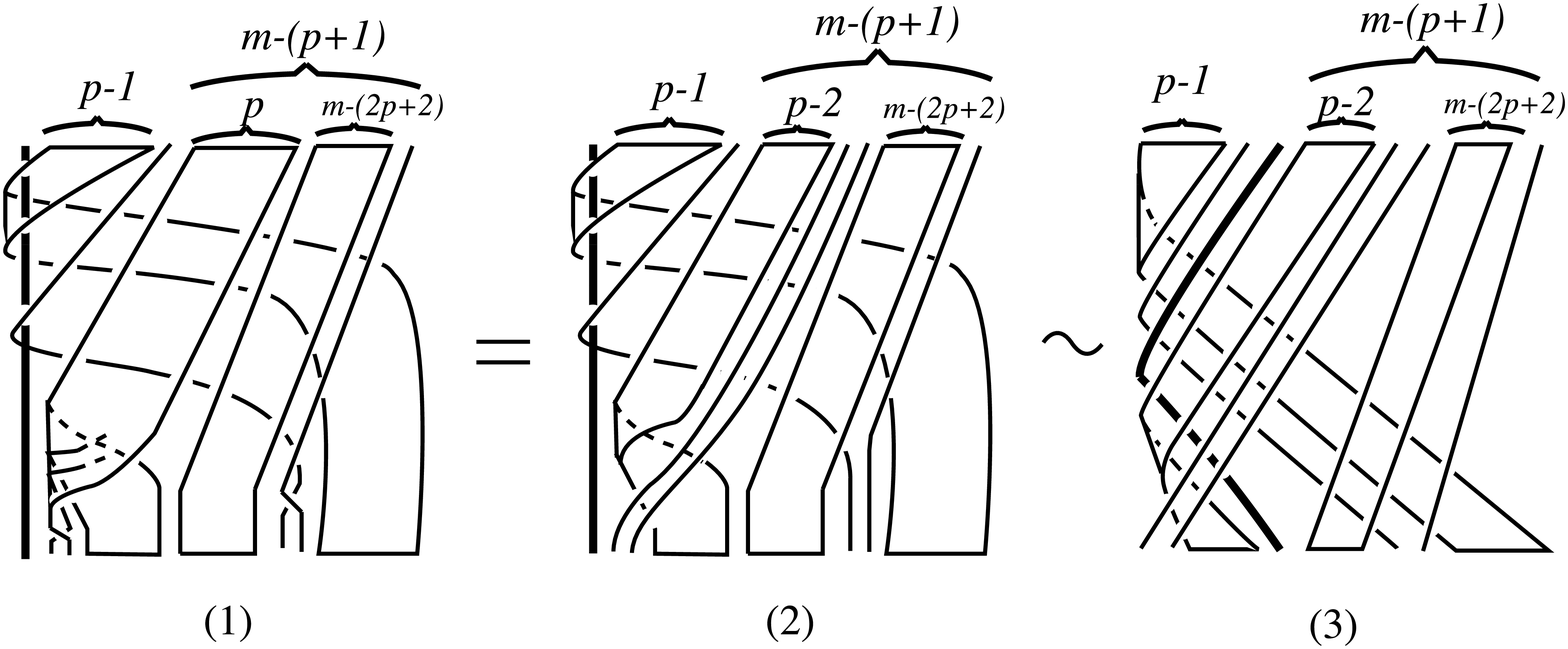}
\caption{conjugate braid  of $T_{m,p}$. }
\label{fig_H-braid-type1}
\end{center}
\end{figure}

First, we show that  the second braid is a horseshoe braid  by using  Figure~\ref{fig_H-braid-type1}.  
This braid is conjugate to the first braid of Figure~\ref{fig_H-braid-type1}  which is equal to the second braid of Figure~\ref{fig_H-braid-type1}. 
 (See the fifth and sixth braid of  Figure~\ref{fig_T_12-4}.)  
The second braid of Figure~\ref{fig_H-braid-type1} is conjugate to the third braid in Figure~\ref{fig_H-braid-type1} which can be embedded in $\mathcal{T}$. 

\begin{figure}[htbp]
\begin{center}
\includegraphics[width=4.5in]{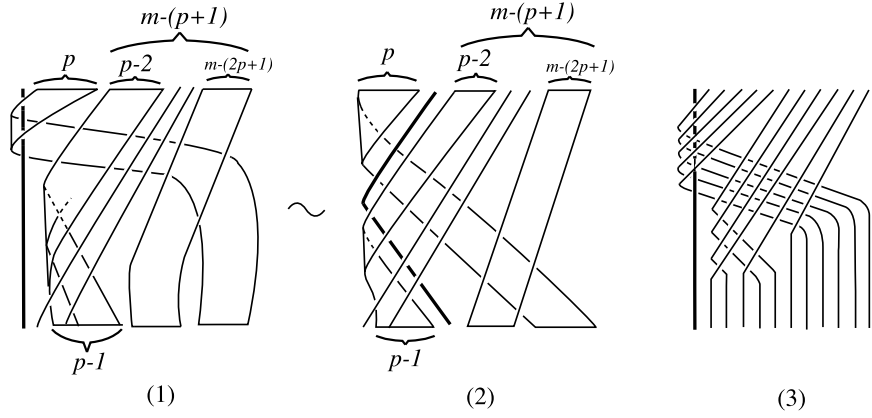}
\caption{(1,2) conjugate braid of $T_{m,p}$. (3) conjugate braid of $T_{12,5}$.}
\label{fig_H-braid-type2}
\end{center}
\end{figure}

Second, we show the third braid of Figure~\ref{fig_horse-braid} is a horseshoe braid by using  Figure~\ref{fig_H-braid-type2}. 
This braid is conjugate to the first braid of Figure~\ref{fig_H-braid-type2}. 
(For example, $T_{12,5}$ is conjugate to the third braid of Figure~\ref{fig_H-braid-type2}.) 
It is easy to see that the first braid of Figure~\ref{fig_H-braid-type2} 
 is conjugate to the second braid of Figure~\ref{fig_H-braid-type2} which can be embedded in $\mathcal{T}$. 
$\Box$

\subsection{Alternative proof of Theorem~\ref{thm_asymp}(2)}
\label{subsection_another-proof}

In this section, we give an alternative proof of Theorem~\ref{thm_asymp}(2).  
\medskip

\begin{figure}[htbp]
\begin{center}
\includegraphics[width=3.8in]{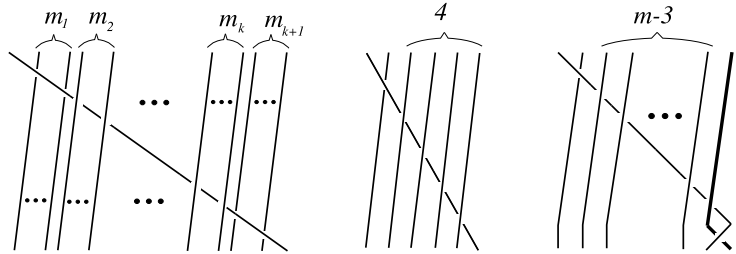}
\caption{(left) $\beta_{(m_1,m_2, \cdots, m_{k+1})}$. (center) $\beta_{(1,4)}$. (right) $\widehat{\beta}_{(1,m-3)}$.}
\label{fig_general_braid}
\end{center}
\end{figure}

\begin{figure}[htbp]
\begin{center}
\includegraphics[width=5in]{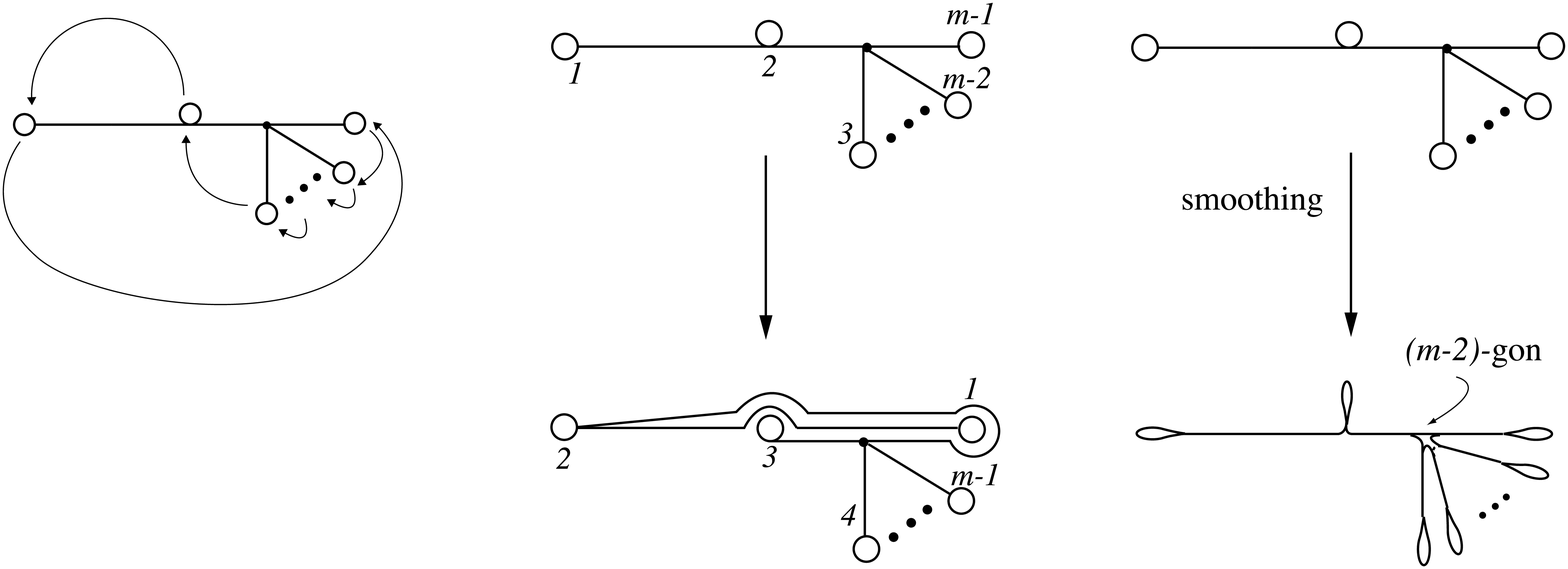}
\caption{(left) transition of peripheral edge. (center) graph map. (right) train track.}
\label{fig_log2}
\end{center}
\end{figure}

\noindent
{\it Proof of Theorem~\ref{thm_asymp}(2).} 
By Proposition~\ref{prop_computation-p}, we have seen that 
$T_{m,1}$ represents the monodromy of $(m-2) \alpha+ \beta \in C_{\Delta_1}({\Bbb Z})$. 
For the proof, it is enough to show that 
\begin{equation}
\label{equation_enough}
\lim_{m (\in {\Bbb N}) \to \infty} \mathrm{ent}((m-2) \alpha+ \beta)= \log2.
\end{equation}
The reason is as follows. 
The equality in (\ref{equation_enough}) implies that 
$$\displaystyle\lim_{x \to \infty} \mathrm{ent}(x \alpha+ \beta)= \log2$$ 
by the continuity of  $\mathrm{ent}(\cdot)$. 
Therefore 
$$\mathrm{ent}(x \alpha+ y \beta) = \tfrac{1}{y} \mathrm{ent}(\tfrac{x}{y} \alpha+ \beta ) \to \tfrac{\log 2}{y}\hspace{2mm}
\mbox{as\ } x \to \infty.$$
Now we show  (\ref{equation_enough}). 
Let $\beta_{(m_1,m_2, \cdots, m_{k+1})}$ be a family of braids depicted  in Figure~\ref{fig_general_braid} 
for each integer $k \ge 1$ and each integer $m_i \ge 1$. 
By \cite[Theorem~1.2]{KT}, 
these braids are all pseudo-Anosov and 
 the dilatation of  $\beta_{(m_1,  \cdots, m_{k+1})}$ is the largest real root of the Salem-Boyd polynomial 
$$ t^{m_{k+1}} R_{(m_1, \cdots, m_k)}(t)+ (-1)^{k+1} {R_{(m_1, \cdots, m_k)}}_*(t),$$
where $R_{(m_1, \cdots, m_i)}(t)$ is given inductively as follows: 
For $2 \le i \le k$, 
\begin{eqnarray*}
R_{(m_1)}(t) &=& t^{m_1+1} (t-1) -2t, 
\\
R_{(m_1, \cdots, m_i)}(t)&=& t^{m_i} (t-1) R_{(m_1, \cdots, m_{i-1})}(t)  + (-1)^i 2t {R_{(m_1, \cdots, m_{i-1})}}_*(t). 
\end{eqnarray*}
In particular, the dilatation of $\beta_{(1,m-3)}= \sigma_1^{-1}\sigma_2 \sigma_3 \cdots \sigma_{m-2}  \in B_{m-1}$ is the largest root of 
\begin{equation}
\label{equation_dil-poly}
t^{m-2} R_{(1)}(t) + (R_{(1)})_*(t),
\end{equation}
where $ R_{(1)}(t)= t(t+1)(t-2)$. 
By Lemma~\ref{lem_asymptotic-root1}, 
the dilatation of $\beta_{(1,m-3)} $ converges to $2$ as $m \to \infty$. 
The polynomial (\ref{equation_dil-poly}) comes from the graph map shown in Figure~\ref{fig_log2}(center). 
This is the induced graph map for  $\beta_{(1,m-3)} \in B_{m-1}$. 
The polynomial   (\ref{equation_dil-poly}) is the characteristic polynomial of the transition matrix for the graph map. 
The smoothing of the graph gives rise to the train track associated to $\beta_{(1,m-3)}$ (Figure~\ref{fig_log2}(right)). 
Since the train track contains an $(m-2)$-gon, a pseudo-Anosov homeomorphism $\Phi_{\beta_{(1,m-3)} }$ which represents the mapping class $\beta_{(1,m-3)} $ 
has an $(m-2)$-pronged singularity, say $p$, in the interior of the punctured disk. 
By puncturing the point $p$, one obtains a pseudo-Anosov homeomorphism $\widehat{\Phi}_{\beta_{(1,m-3)} }$. 
It is easy to see that the mapping class $[\widehat{\Phi}_{\beta_{(1,m-3)} }]$ is given by 
$$\widehat{\beta}_{(1,m-3)} = \sigma_1^{-1}\sigma_2 \sigma_3 \cdots \sigma_{m-2} \sigma_{m-1}^2 \in B_m$$
with the same dilatation as $\beta_{(1,m-3)}$. 
Since $\widehat{\beta}_{(1,m-3)}$ is conjugate to the braid $T_{m,1}$, the dilatation $\lambda(T_{m,1})$ converges to $2$ as $m$ goes to $\infty$. 
This completes the proof. 
$\Box$

\noindent
Department of Mathematical and Computing Sciences, Tokyo Institute of Technology, 
  Ohokayama, Meguro  Tokyo 152-8552 Japan \\
 E-mail address: kin@is.titech.ac.jp
 \medskip
 
 \noindent
 Department of Mathematical and Computing Sciences,  Tokyo Institute of Technology, 
  Ohokayama, Meguro  Tokyo 152-8552 Japan \\
 E-mail address: takasawa@is.titech.ac.jp

\end{document}